\documentclass[graybox]{svmult}

\usepackage{mathptmx}       
\usepackage{helvet}         
\usepackage{courier}        
\usepackage{type1cm}        
\usepackage{makeidx}         
\usepackage{graphicx}        
\usepackage{multicol}        
\usepackage[bottom]{footmisc}

\usepackage{mathrsfs}
\usepackage{amsmath}
\usepackage{amssymb}

\makeindex             

\usepackage{dsfont}   
\usepackage{euscript}
\DeclareMathOperator{\sgn}{sgn}
\newcommand{\lawto}{\xrightarrow{\scriptscriptstyle{\textup{law}}}}
\font\gfont=cmmi10 scaled \magstep{1.5}     
\newcommand{\gdelta}{\hbox{\gfont \char14}}

\begin{document}

\title*{Fluctuations of observables in dynamical systems: from limit theorems to concentration inequalities}
\titlerunning{From limit theorems to concentration inequalities}
\author{Jean-Ren\'e Chazottes}
\institute{Jean-Ren\'e Chazottes \at Centre de Physique Th\'eorique, CNRS-\'Ecole polytechnique, 91128 Palaiseau Cedex, France, \email{jeanrene@cpht.polytechnique.fr}
}
%
%
\maketitle

\abstract*{We start by reviewing recent probabilistic results on ergodic sums in a large class of (non-uniformly) hyperbolic dynamical systems. Namely, we describe the central limit theorem, the almost-sure 
convergence to the Gaussian and other stable laws, and large deviations.\newline
Next, we describe a new branch in the study of probabilistic properties of dynamical systems, namely concentration inequalities.
They allow to describe the fluctuations of very general observables and to get bounds rather than limit laws.\newline
We end up with two sections: one gathering various open problems, notably on random dynamical systems,
coupled map lattices and the so-called nonconventional ergodic averages; and another one giving pointers
to the literature about moderate deviations, almost-sure invariance principle, etc.} 

\abstract{
We start by reviewing recent probabilistic results on ergodic sums in a large class of (non-uniformly) hyperbolic dynamical systems. Namely, we describe the central limit theorem, the almost-sure 
convergence to the Gaussian and other stable laws, and large deviations.\newline
Next, we describe a new branch in the study of probabilistic properties of dynamical systems, namely concentration
inequalities. They allow to describe the fluctuations of very general observables and to get bounds rather than limit laws.\newline
We end up with two sections: one gathering various open problems, notably on random dynamical systems,
coupled map lattices and so-called nonconventional ergodic averages; and another one giving pointers to the literature about moderate deviations, almost-sure invariance principle, etc.
}

\bigskip

\bigskip

\setcounter{minitocdepth}{1}
\dominitoc

\newpage

\section{Introduction}

The aim of the present chapter is to roughly describe the current state of the theory of {\em statistical} or {\em probabilistic}
properties of `chaotic' dynamical systems. We shall restrict ourselves to {\em discrete-time} dynamical systems, although many of the results we review have their counterparts in flows. The basic setting is thus a
state space $\Omega$ (typically a piece of $\mathds{R}^d$) and a map $T:\Omega\circlearrowleft$. The orbit of an
initial condition $x_0$ is the sequence of points $x_0,x_1=Tx_0,x_2=Tx_1,\ldots$
or $\{T^k x_0; k=0,1,\ldots\}$
(where $T^k$ is the $k$-fold composition of $T$ with itself).

The core of the probabilistic  approach is the description of asymptotic {\em time-averages} of `observables', that is, functions $f:\Omega\to\mathds{R}$. 
This implies that transients become irrelevant, although transient effects may cause formidable problems in practice.
The corner stone of this approach is Birkhoff's ergodic theorem. 
It tells us that, {\em given} a measure $\mu$ left invariant by $T$, `the asymptotic {\em time-average} of $f$ coincide with
the space-average $\int f \textup{d}\mu$', except on a set of measure zero with respect to this measure. The drawback of this result is that chaotic systems typically possess uncountably many invariant ergodic measures. Is there a `natural' choice ?

In this chapter, we focus on {\em dissipative} systems whose orbits settle on an attractor which has typically a volume (Lebesgue measure) equal to zero. In these systems, the dynamics contracts volumes but generally not in all directions: some directions may
be stretched, provided some others are so much contracted that the final volume is smaller than the initial volume.
This implies that, even in a dissipative system, the motion after transients may be unstable within the attractor.
This instability manifests itself by an exponential separation of orbits, as time goes on, of points which initially
are very close to each other on the attractor. The exponential separation takes place in the direction of stretching.
Such an attractor is called {\em chaotic}. Of course, since the attractor is bounded, exponential separation can only
hold as long as distances are small.

A famous attractor is the H\'enon attractor generated by a two-dimensional map with two parameters. For some parameters,
it is easy to numerically produce a `picture' of the attractor. The standard way to make it is to pick `at random'
an initial condition in the basin of the attractor and to plot the first thousand iterates of its orbit (see Fig. \ref{henon-fig}).
On the one hand, why does what is observed has something to do with the attractor since, as noticed above,
it has zero volume ? On the other hand, we know that orbits of the H\'enon map are not all the same: some are periodic,
others are not; some come closer to the `turns' than others. We also know from experience that (for a fixed $T$)
one gets essentially the same picture independent of the choice of initial condition. Is there a mathematical explanation for
this ?

These questions motivated the idea of {\em Sinai-Ruelle-Bowen} or SRB measures. Our computer picture can be thought
as the picture of a probability measure giving mass $1/n$ to each point in an orbit of length $n$. Let $\delta_x$
be the point mass at $x$. Is there a (probability) measure $\mu$ with the property that $\frac{1}{n} \sum_{i=0}^{n-1}
\gdelta_{T^i(x)}\to \mu$ for `most' choices of initial conditions $x$, that would explain why our pictures look similar ?
If such a measure does exist, it has very special properties: like all invariant probability measures, it must be supported
on the attractor, but it has the peculiar ability to influence orbits starting from various parts of the basin, including points
rather far away from the support of $\mu$. In some sense, SRB measures are the observable or physical measures.

Mathematically speaking, the theory of chaotic attractors began with the ergodic theory of {\em differentiable dynamical
systems}, more specifically the theory of hyperbolic dynamical systems, where geometry plays a prominent role.
The first systems studied in the 1960-70's were the so-called Anosov and Axiom A systems
which are `uniformly' hyperbolic and in some sense the most chaotic systems.
The main results were obtained by Sinai, Ruelle and Bowen. They essentially relied on the fact that, for such
systems, it is possible to construct Markov partitions enabling one to identify points in the state space with configurations in one-dimensional lattice systems of statistical mechanics \cite{bowen}.\newline
The 1970's brought new outlooks and new challenges. With the aid of computer graphics, an abundance of examples
showed up whose dynamics is dominated by expansions and contractions, but which do not satisfy the stringent requirements
of Axiom A systems. H\'enon's attractor mentioned above is a typical example. This led to a more comprehensive theory
dealing with {\em non-uniformly} hyperbolic dynamical systems developed abstractly by Pesin and others \cite[Chap. 2]{handbookb}. A breakthrough was made by L.-S. Young at the end of the 1990's \cite{young1,young2}. She proposed a more `phenomenological' approach to describe in a unified framework many examples of systems with a `localized' source of non-hyperbolicity.  In particular, this provided tools to prove the existence of an SRB measure for the H\'enon attractor (for a set of parameters with positive measure), see \cite{BY}. In this chapter, we shall focus on the class of systems defined by Young.

Once we know that our dynamical system $(\Omega,T)$ admits an SRB measure, we can ask for its probabilistic
properties. Indeed, it can be viewed as a stationary stochastic process: the
orbits $(x,Tx,\ldots)$, where $x$ is distributed according to $\mu$, generate a stationary process whose finite-dimensional marginals are the measures $\mu_n$  on $\Omega^n$ given by 
\[
\textup{d}\mu_n(x_0,\ldots,x_{n-1})=\textup{d}\mu(x_0)\delta_{x_1=Tx_0}\cdots \delta_{x_{n-1}=Tx_{n-2}}.
\]
This is not a product measure but the idea is that, if the system is chaotic enough, $T^k x$ is more or less
independent of $x$ provided $k$ is large, making the process $(x,Tx,\ldots)$ behave like an independent process.\newline
Given any observable $f:\Omega\to\mathds{R}$,  one can generate a process $\{X_n=f\circ T^n; n\geq 0\}$  on
the probability space $(\Omega,\mu)$.  The ergodic sum $S_n f(x)=f(x)+f(Tx)+\cdots+f(T^{n-1}x)$
is thus the partial sum of the process $\{X_n; n\geq 0\}$ and one can ask various natural questions.
For instance, what is the typical size
of fluctuations of $\frac{1}{n}S_n f(x)$ around $\int f\textup{d}\mu$ ? What is the probability that
$\frac{1}{n}S_n f(x)$ deviates from $\int f\textup{d}\mu$ by more than some prescribed value ?
Does $S_nf$, appropriately renormalized, converge
in law ? In other words, can we prove a central limit theorem ? Can we get a description of large deviations ?
Can we have Gaussian but also non-Gaussian limit laws ? This kind of results are called {\em limit theorems}.

There are many quantities describing a dynamical system which can be in principle computed by observing
its orbits. But the corresponding estimators are not as simple as ergodic sums of suitably chosen observables.
A prominent example (see below for details) is the periodogram which is related to the power spectrum.
Therefore it is desirable to have a tool which allows to quantify fluctuations of fairly general observables for
finite-length orbits. This is the scope of concentration inequalities, a new branch in the study of probabilistic
theory of dynamical systems (and a quite recent branch of Probability theory as well \cite{ledoux}). The aim
of concentration inequalities is to quantify the size of the deviations of an observable $K(x,Tx,\ldots,T^{n-1}x)$
around its expectation, where $K:\Omega^n\to\mathds{R}$ is an observable of $n$ variables of an arbitrary expression.
An ergodic sum is a very special case of such an observable and we shall see below various examples.
What is imposed on $K$ is sufficient smoothness (Lipschitz property).
Depending on the `degree of chaos' in the system, the deviations of $K$ with respect to its expectation
can have an extremely small probability.

From the technical viewpoint, the tool of paramount importance is the {\em transfer} or {\em Ruelle's Perron-Frobenius
operator}. This is the spectral approach to dynamical systems. We refer to book of Baladi \cite{viviane} and to the
lecture notes of Hennion and Herv\'e \cite{HH} for a throughout exposition. 

Our purpose is to give a sample of recent results on the fluctuations of observables in
the ergodic theory of non-uniformly hyperbolic dynamical systems.
Needless to say that the overwhelming list of works in this area renders
futile any attempt at an exhaustive or even comprehensive treatment within the confines 
of this chapter. Hopefully, this chapter provides a panoramic view of this subject.
We also provide a list of directions for further research.

Before describing the contents of this chapter, a few words are in order about
the bibliography. We urge the reader to consult \cite{compil} in which are gathered
landmark papers illustrating the history and development of the notions of chaotic
attractors and their `natural' invariant measures. For numerical implementations
of the theory, it is still worth reading the review paper by Eckmann and Ruelle \cite{ER}.
A more recent reference, dealing both with theoretical and numerical aspects is the
book by Collet and Eckmann \cite{pierre-jp}. Needless to say that the potential
list of references is gigantic. Limitation of space and time forced us painfully
to exclude many relevant papers. As a matter of principle, and whenever possible,
we refer to the most recent articles which contain relevant pointers to the literature.
We apologize for omissions.

\noindent{\em Layout of the chapter}.
In Section \ref{sec:1} we describe the probability approach to dynamical systems and recall Birkhoff's
ergodic theorem. In Section \ref{sec:2} we describe the class of hyperbolic dynamical systems we
will be working with. In particular, we quickly describe Young towers and SRB measures, and give
several examples which will be used throughout the chapter. Section \ref{sec:3} is devoted to
mixing (decay of correlations) and limit theorems, namely: the central limit theorem, convergence to non-Gaussian laws,
exponential and sub-exponential large deviations, and convergence in law made almost sure.
Section \ref{sec:4} is concerned with
concentration inequalities and some of their applications. In Section \ref{sec:5} we provide a list of open
problems and questions related to random dynamical systems, coupled map lattices, partially hyperbolic systems,
and the Erd\"os-R\'enyi law. We end with a section where we quickly survey results not detailed in the main
text. This includes Berry-Esseen theorem, moderate deviations and the almost-sure invariance principle.

\section{Generalities}
\label{sec:1}

We state some general definitions and recall Birkhoff's ergodic theorem.

\subsection{Dynamical systems and observables}

In this chapter, by `dynamical system' we mean a
deterministic dynamical system with discrete time, that is,
a transformation $T:\Omega\circlearrowleft$ of its state
space (or phase space)  $\Omega$ into itself. For the sake of concreteness, one
can think of $\Omega$ as a compact subset of $\mathds{R}^d$. Mathematically speaking,
one can deal with a compact riemannian manifold.

Every point $x\in\Omega$ represents a possible state of the system. If the system is in state $x$, then
it will be in state $T(x)$ in the next moment of time. Given the current state
$x=x_0\in\Omega$, the sequence of states
\[
x_1=T x_0,\; x_2=T x_1,\ldots,\; x_n=Tx_{n-1},\ldots
\]
represents the entire future or forward orbit of $x_0$. We have $x_n=T^n x_0$, where
$T^n$ is the $n$-fold composition of $T$ with itself. If the map $T$ is invertible, then
the past of $x_0$ can be determined as well ($x_{-n}=T^{-n}x_0$).

In applications, the actual states $x_n\in\Omega$ are often not observable. Instead,
we usually observe the values $f(x_n)$ taken by a function $f$ on $\Omega$, usually called
an observable. One can be thought of $f$ as an instantaneous measurement of the system. For the
sake of simplicity, we consider $f$ to be real-valued. 

More  generally, we may observe the system from time $0$ up to time $n-1$ and 
associate to $x,Tx,\ldots,T^{n-1}x$ a real number $K(x,Tx,\ldots,T^{n-1}x)$. In the language of statistics,
$K:\Omega^n\to\mathds{R}$ is called an estimator.
The fundamental example is the Ces\`aro or ergodic average of an `instantaneous' observable
$f:\Omega\to\mathds{R}$ along an orbit up to time $n-1$:
$K_0(x,Tx,\ldots,T^{n-1}x):=(f(x)+f(Tx)+\cdots+f(T^{n-1}x))/n$. This is an example of an additive observable. There are
many natural examples which are not as simple.
An important example is the periodogram used to estimate the power spectrum of
a `signal' $\{f(x_k); k=0,\ldots,n-1\}$. We give its definition below as well as other examples; see section \ref{subsec:applications}.

\subsection{Dynamical systems as stochastic processes}

Ergodic theory is concerned with measure-preserving transformations, meaning that the map $T$ preserves
a probability measure $\mu$ on $\Omega$: for any measurable subset $A\subset \Omega$ one has
$\mu(A)=\mu(T^{-1}(A))$, where $T^{-1}(A)$ denotes the set of points mapped into $A$. The invariant
measure $\mu$ describes the distribution of the sequence $\{x_n=T^{n-1}(x_0)\}$ for typical initial states
$x_0$. This vague statement is made precise by Birkhoff's ergodic theorem; see below. For a large
class of non-uniformly hyperbolic systems, there is a `natural' invariant measure, the so-called
Sinai-Ruelle-Bowen measure (SRB measure for short).

A measure-preserving dynamical system is thus a probability space
$(\Omega,\mathscr{B},\mu)$ endowed with a transformation $T:\Omega\circlearrowleft$ 
leaving $\mu$ invariant.  An important notion is that of an ergodic dynamical system.
The invariant measure $\mu$ is said to be ergodic (with respect to $T$) whenever 
$T^{-1}(E)=E$ implies $\mu(E)=0$ or $\mu(E)=1$. Equivalently, ergodicity means
that any invariant function $g:\Omega\to\mathds{R}$ is $\mu$-almost everywhere constant.
That $g$ be invariant means that $g=g\circ T$. In the measure-theoretic sense, ergodic measures are
indecomposable and any invariant measure can be disintegrated into its ergodic components \cite{KH}.

A measure-preserving dynamical system can be viewed as a stochastic process: the
orbits $(x,Tx,\ldots)$, where $x$ is distributed according to $\mu$, generate a stationary process whose finite-dimensional marginals
are the measures $\mu_n$  on $\Omega^n$ given by 
\[
\textup{d}\mu_n(x_0,\ldots,x_{n-1})=\textup{d}\mu(x_0)\delta_{x_1=Tx_0}\cdots \delta_{x_{n-1}=Tx_{n-2}}.
\]
This is not a product measure but the idea is that, if the system is chaotic enough, $T^k x$ is more or less
independent of $x$ provided $k$ is large, making the process $(x,Tx,\ldots)$ behave like an independent process.

Given an observable $f:\Omega\to\mathds{R}$, $X_k=f\circ T^k$, for each $k\geq 0$, is a random variable on
the probability space $(\Omega,\mathscr{B},\mu)$. The family $\{X_n; n\geq 0\}$ is a real-valued stationary process.
The ergodic sum $S_n f(x)=f(x)+f(Tx)+\cdots+f(T^{n-1}x)$ is thus the partial sum of the process $\{X_n; n\geq 0\}$.

We shall make no attempt to define precisely what a chaotic dynamical system is. From the point of view
of this chapter, we can vaguely state that it is a system such that, for sufficiently nice observables $f$, the process
$\{f\circ T^k\}$ behave as an i.i.d.\footnote{i.i.d. stands for `independent and identically distributed'} process. Along the way, this crude statement will be refined.

\subsection{Birkhoff's Ergodic Theorem}\label{subsec:birk}

The fundamental theorem in ergodic theory is Birkhoff's ergodic theorem which is a far reaching generalization of Kolmogorov's strong law of large numbers for an independent process \cite{krengel}.
\begin{svgraybox}
\begin{theorem}[Birkhoff's ergodic theorem]\label{birkhoff}
\leavevmode\\
Let $(\Omega,\mathscr{B},\mu)$ be a dynamical system and $f:\Omega\to\mathds{R}$ be an integrable
observable \textup{(}$\int |f|\textup{d}\mu<\infty$\textup{)}. Then
\[
\lim_{n\to\infty} \frac{1}{n}S_n f(x)=f^*(x),\;\;\mu-\textup{almost surely and in}\;L^1(\mu),
\]
where the function $f^*$ is invariant ($f^*=f^*\circ T$, $\mu$-a.s.) and such that
$\int f^* \textup{d}\mu=\int f\textup{d}\mu$.\newline
\item[]If the dynamical system is ergodic, then $f^*$ is $\mu$-almost surely a constant, whence
\[
\lim_{n\to\infty} \frac{1}{n}S_n f(x)=\int f\textup{d}\mu,\;\;\mu-\textup{almost surely}.
\]
\end{theorem}
\end{svgraybox}

\begin{remark}
The previous theorem, spelled out for an integrable stationary ergodic process $\{X_n\}$, reads
$n^{-1}\sum_{j=0}^{n-1} X_j\xrightarrow{} \mathds{E}[X_0]$ almost surely.
In the non-ergodic case convergence is to the conditional expectation of $X_0$
with respect to the $\sigma$-algebra of invariant sets, see \cite{krengel} for details.
\end{remark}

Very often, $\Omega$ is compact and it is not difficult to show that there exists a measurable set of $\mu$-measure one such that, in the ergodic case, 
$g:\Omega\to\mathds{R}$
\[
\lim_{n\to\infty} \frac{1}{n}S_n g(x)=\int g\textup{d}\mu
\]
for any continuous observables.
Equivalently, this means that the \emph{empirical measure} of $\mu$-almost every $x$ converges towards $\mu$ in the
vague (or weak-$^*$) topology:
\[
\frac{1}{n}\sum_{j=0}^{n-1} \delta_{T^j x} \xrightarrow{\scriptscriptstyle{\textup{vaguely}}}\mu\quad\textup{almost surely}.
\]
The advantage of Birkhoff's ergodic theorem is its generality. Its drawback is that a chaotic system has in general
uncountably many distinct ergodic measures. Which one do we choose ? We shall see later on that the idea of a
Sinai-Ruelle-Bowen measure provides an answer. 

\subsection{Speed of convergence and fluctuations}

It is well known that not much can be said about the speed of convergence of the ergodic average to its limit in
Theorem \ref{birkhoff}. First of all, one cannot know in practice if we are observing a typical orbit for which convergence
indeed occurs. But even if we knew that we have a typical orbit, it can be shown that the convergence can be arbitrarily
slow (see for instance \cite{kachu} for a survey).\newline
To obtain more informations about the fluctuations of ergodic sums around their limit, we need a probabilistic formulation.
Maybe the most natural question is the following: 

\medskip
\begin{minipage}{10cm}
what is the speed of convergence to zero of the probability that
the ergodic average differs from its limit by more than a prescribed value ? 
\end{minipage}
\medskip

\noindent Formally, we want to know the speed of convergence to zero of
\[
\mu\left\{x : \Big|\frac{1}{n}S_nf(x)-\int f\textup{d}\mu\Big|>t \right\}
\]
for $t>0$ small enough and for a large class of continuous observables $f$. (By Birkhoff's ergodic
theorem, all what we know is that this probability goes to $0$ as $n$ goes to infinity.)
In probabilistic terminology, we want to know the speed of convergence \emph{in probability} of
ergodic averages to their limit. By analogy with bounded i.i.d. processes, this speed should be
exponential for `sufficiently chaotic' systems. We shall see that it can be only polynomial
when mixing is not strong enough.

Another natural issue is to determine the order of \emph{typical values} of $S_n f-n\int f\textup{d}\mu$.
By analogy with a square-integrable  i.i.d. process, one can expect this order to be $\sqrt{n}$, and, more precisely,
 that a central limit theorem may hold. We shall see that this is indeed the case for `nice observables'
and sufficiently chaotic systems. When chaos is `too weak', the central limit theorem may fail and
the asymptotic distribution may be non-Gaussian.

The previous issues are formulated in terms of limit theorems and concern ergodic sums. From the
point of view of applications, an important problem is to estimate the probability of deviation
of a general observable $K(x,Tx,\ldots,T^{n-1}x)$ from its expected value. Formally, we ask
if it is possible to find a positive function $b(n,t)$ such that
\[
\mu\left\{x\in\Omega : \Big|K(x,Tx,\ldots,T^{n-1}x)-\int K(y,Ty,\ldots,T^{n-1}y)\textup{d}\mu(y)\Big|>t \right\}\leq b(n,t)
\]
for any $t>0$ and for any $n\in\mathds{N}$, with $b(n,t)$ depending on $K$. When $b(n,t)$ decreases `rapidly' with $t$ and $n$,
this means that $K(x,Tx,\ldots,T^{n-1}x)$ is `concentrated' around its expected value. It turns out that
when the dynamical system is `chaotic enough', this concentration phenomenon is very sharp.

To be able to answer the kind of previous questions, we shall need to make hypotheses on the dynamical systems
as well as on the class of observables. Usually, H\"older continuous functions are suitable.

\section{Dynamical systems with some hyperbolicity}
\label{sec:2}

We quickly and roughly describe the class of dynamical systems for which one can prove various probabilistic
results. These systems are used to model deterministic chaos which is caused by dynamic instability,
or sensitive dependence on initial conditions, together with the fact that orbits are confined in a compact
region.

\subsection{Hyperbolic dynamical systems}

The basic model for sensitive dependence on initial conditions is that of a uniformly expanding map $T$
on a riemannian compact manifold $\Omega$: $T$ is smooth and there are constants $C>0$
and $\lambda>1$ such that for any $x\in\Omega$ and $v$ in the tangent space at $x$ and for any $n\in\mathds{N}$
\[
\| DT^n(x) v\|\geq C \lambda^n \|v\|.
\]
The prototypical example is $T(x)=2x$ (mod $1$) on $\Omega=S^1$ (the unit circle), which is usually
identified with the interval $[0,1)$. The Lebesgue measure is invariant in this case.

Uniformly hyperbolic maps have the property that at each point $x$ the tangent space is a direct sum of two
subspaces $E^u_x$ and $E^s_x$, one of which is expanded ($\| DT^n(x) v\|\geq C \lambda^n \|v\|$ for $v\in E^u_x$) and the other contracted ($\| DT^n(x) v\|\leq C \lambda^{-n} \|v\|$ for $v\in E^s_x$). The prototypical example is Arnold's cat map $(x,y)\mapsto (2x+y,x+y)$ (mod $1$) of the unit torus.

Non-uniform hyperbolicity refers to the fact that $C=C(x)>0$ and $\lambda=\lambda(x)>1$
almost-everywhere: in words, the constants depend on $x$ and they have nice properties only
on a set a full measure. For instance, the presence of a single point where $\lambda(x)=1$ already causes
important difficulties (the fundamental example being an interval map with an indifferent fixed point at $0$).
Another instance of loose of uniform hyperbolicity is when there is a point where the differential
of $T$ vanishes ({\em e.g.}, the quadratic map or the H\'enon map). A third typical situation is when the differential has discontinuities. This is the case for the Lozi map and billiards, for instance.

\subsection{Attractors}

We are especially interested in dissipative
systems with an attractor, that is, volume-contracting maps $T$ with an attractor 
$\Lambda$. By an attractor we refer to a compact invariant set with the property
that all points in a neighborhood $U$ of $\Lambda$ (called its basin) are
attracted to $\Lambda$ ({\em i.e.} for any $x\in U$, $T^n x\to \Lambda$ as $n\to\infty$).

The prototype of a hyperbolic attractor is an Axiom A attractor. It is a smooth map $T$ with an attractor $\Lambda$ on which $T$
is uniformly hyperbolic. These systems can be viewed as subshifts of finite type
by using a Markov partition: one can assign to each point a bi-infinite symbol
sequence describing its itinerary. This sequence can be thought of as a configuration
in a one-dimensional statistical mechanical system. Special measures, called
SRB measures (see next section) can be constructed by pulling back adequate Gibbs measures which are invariant by the shift map; see \cite{bowen} and \cite[Chap. 4]{handbooka}.

H\'enon's attractor is a genuinely non-uniformly hyperbolic attractor which resisted to mathematical analysis till
the 1990's.

\subsection{Sinai-Ruelle-Bowen measures}

We shall not define precisely Sinai-Ruelle-Bowen (SRB for short) measures but content ourselves by saying that they
are the invariant measures most compatible with volume (Lebesgue measure) when volume is not preserved.
Technically speaking, they have absolutely continuous conditional measures along unstable manifolds and
a positive Lyapunov exponent.
They provide a mechanism for explaining how local instability on attractors can produce coherent
statistics for orbits starting from large sets in the basin. In particular, an SRB measure $\mu$ is `observable'
in the following sense: there exists a subset $V$ of the basin of attraction with positive Lebesgue measure such that for
any continuous observable $f$ on $\Omega$ and any initial state $x\in V$ we have
\[
\lim_{n\to\infty} \frac{1}{n} \sum_{j=0}^{n-1} f(T^j x)=\int f\textup{d}\mu,
\]
or, more compactly
\[
\frac{1}{n}\sum_{j=0}^{n-1} \delta_{T^j x} \xrightarrow{\scriptscriptstyle{\textup{vaguely}}}\mu.
\]
The point of this property is that the set of `good states' has positive Lebesgue measure
although the measure $\mu$ is concentrated on the attractor which has zero Lebesgue 
measure. (Notice that this property does not follow from Birkhoff's ergodic theorem.) 

For one-dimensional maps, absolutely continuous invariant measures (with respect to Lebesgue
measure) are examples of SRB measures.

Roughly speaking, the approach to non-uniformly hyperbolic systems of L.-S. Young, which will be sketched below,
can be considered as `phenomenological' in the sense that it aims at modeling concrete dynamical behaviors
observed in various examples. An `axiomatic approach' can be followed which
seeks to relax the conditions that define Axiom A systems in the hope of systematically
enlarging the set of maps with SRB measures. For an account on this second approach,
we refer to \cite[Chap. 2]{handbookb}.
For a nice and non-technical survey on SRB measures, we recommend reading \cite{young3}.

\subsection{Dynamical systems modeled by a Young tower}

In the 1970s, many examples were numerically observed whose dynamics
are dominated by expansions and contractions but which do not meet
the stringent requirements of Axiom A systems. The most famous example
is likely the H\'enon mapping which displays a `strange attractor' for certain parameters. Such examples remained mathematically intractable until the 1990's. 

L.-S. Young developped a general scheme to study the probabilistic properties of a class of `predominantly hyperbolic' dynamical
systems, including the H\'enon attractor and other famous examples.
Very roughly the picture is as follows.
The general set up is that $T:\Omega\circlearrowleft$ is a nonuniformly hyperbolic system in the sense of
Young \cite{young1,young2} with a return time function $R$ that decays either exponentially \cite{young1}, or polynomially \cite{young2}.
In particular, $T:\Omega\circlearrowleft$  is modeled by a Young tower constructed over a `uniformly hyperbolic' base $Y\subset \Omega$. The degree of non-uniformity is measured by the return time function $R:Y\to \mathds{Z}^+$ to the base.

\begin{figure}
\sidecaption[t]
\includegraphics[scale=.29]{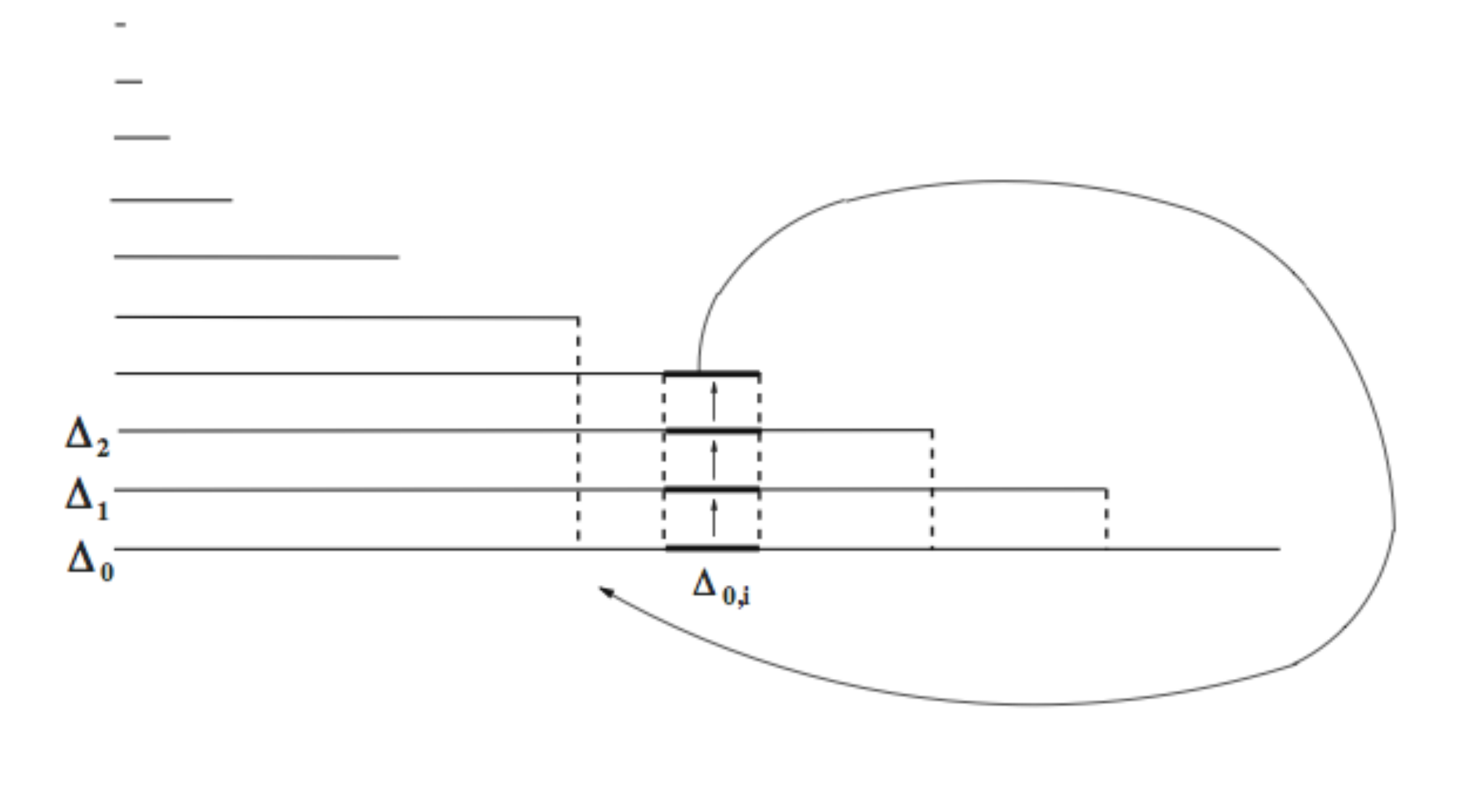}
\caption{Schematic representation of the tower map $F:\Delta\circlearrowleft$.}
\label{dessintour}
\end{figure}

More precisely, by a classical construction in ergodic theory, one can construct from $(Y,T^R)$ an
extension $(\Delta,F)$, called a Young tower in the present setting.
In particular, there exists a continuous map $\pi:\Delta\to\Omega$ such that
$\pi\circ F= T\circ \pi$. In general $\pi$ need not be one-to-one or onto.
One can visualize a tower by writing that $\Delta=\cup_{\ell=0}^\infty \Delta_\ell$ where $\Delta_\ell$ can be identified
with the set $\{x\in Y : R(x)>\ell\}$, that is, the $\ell$-th floor of the tower. In particular, $\Delta_0$
is identified with $Y$. The dynamics in the tower is as follows: each point $x\in\Delta_0$
moves up the tower until it reaches the top level above $x$, after which it returns to
$\Delta_0$, see Fig. \ref{dessintour}. Moreover, $F$ has a countable Markov partition
$\{\Delta_{\ell,j}\}$ with the property that $\pi$ maps each $\Delta_{\ell,j}$ injectively
onto $Y$, which has a hyperbolic product structure. Each of the local unstable manifolds defining the
product structure of $\pi(\Delta_0)$ meet $\pi(\Delta_0)$ in a set of positive Lebesgue measure. Further analytic and regularity conditions are imposed. We shall not give further details and refer the reader to \cite{young1,young2} and \cite{jrseb}. 

Systems modeled by Young towers are more flexible than Axiom A systems in that
they are permitted to be non-uniformly hyperbolic: roughly speaking, think of
uniform hyperbolicity as required only for the return map to the base. Reasonable
singularities and discontinuities are also allowed: they do not appear in $Y$.
As we shall see, a number of probabilistic properties of $T:\Omega\circlearrowleft$ are actually captured by the tail
properties of $R$.
The basic result proved in \cite{young1,young2} is the following, where $m^u$ denotes Lebesgue measure
on unstable manifolds.
\begin{theorem}
\begin{svgraybox}
Let $T:\Omega\circlearrowleft$ be a dynamical system modeled by a Young tower.\\
If $\int R \textup{d}m^u<\infty$, then $T$ has an ergodic SRB measure. If $\textup{gcd}\{R_i\}=1$,
there is a unique SRB measure denoted by $\mu$.
\end{svgraybox}
\end{theorem}
Of course, $\int R \textup{d}m^u=\sum_{n\geq 1}m^u\{R>n\}$.
In the sequel, we shall implicitly assume that $\textup{gcd}\{R_i\}=1$, without loss of generality.

\subsection{Some examples}

The best known example of a non-uniformly expanding map of the interval is the
so-called Maneville-Pomeau map modelling intermittency. It is expanding except at $0$ where the slope
of the map is one (neutral fixed point). For the sake of definiteness\footnote{The explicit formula \eqref{MPmap} is not important, what matters is only the local behavior around the fixed point.}, consider the map
\begin{equation}\label{MPmap}
T_\alpha(x)=
\begin{cases}
x+2^\alpha x^{1+\alpha} & \textup{if} \quad x\in[0,1/2)\\
2x-1 & \textup{if}\quad x\in [1/2,1]
\end{cases}
\end{equation}
where $\alpha\in(0,1)$ is a parameter. It is well-known that there is a unique
absolutely continuous invariant measure $\textup{d}\mu(x)=h(x)\textup{d}x$
and $h(x)\sim x^{-\alpha}$ as $x\to 0$. There is a Young tower with base $Y=[1/2,1]$
and $\textup{Leb}\{y\in Y : R(y)>n\}=\mathcal{O}(n^{-1/\alpha})$.

Another fundamental one-dimensional example is given by the quadratic family $T_a:[-1,1]\circlearrowleft$
with $T_a(x)=1-ax^2$, where $a\in[1,2]$, and for which $0$ is a critical point (the slope vanishes). For a set
of parameters of positive Lebesgue measure, this maps preserves a unique absolutely continuous probability
measure. Its density has an inverse square-root singularity. In this example, one can construct a tower map with a return-time function which has an exponentially
decreasing tail.

\begin{figure}
\sidecaption[t]
\includegraphics[scale=.35]{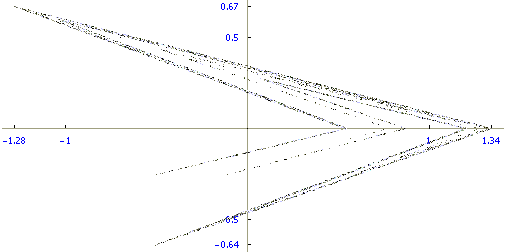}
\caption{Simulation of the Lozi attractor for $a=1,7$ and $b=0,5$.}
\label{A-lozi}
\end{figure}

\begin{figure}[htb!]
\sidecaption[t]
\includegraphics[scale=.09]{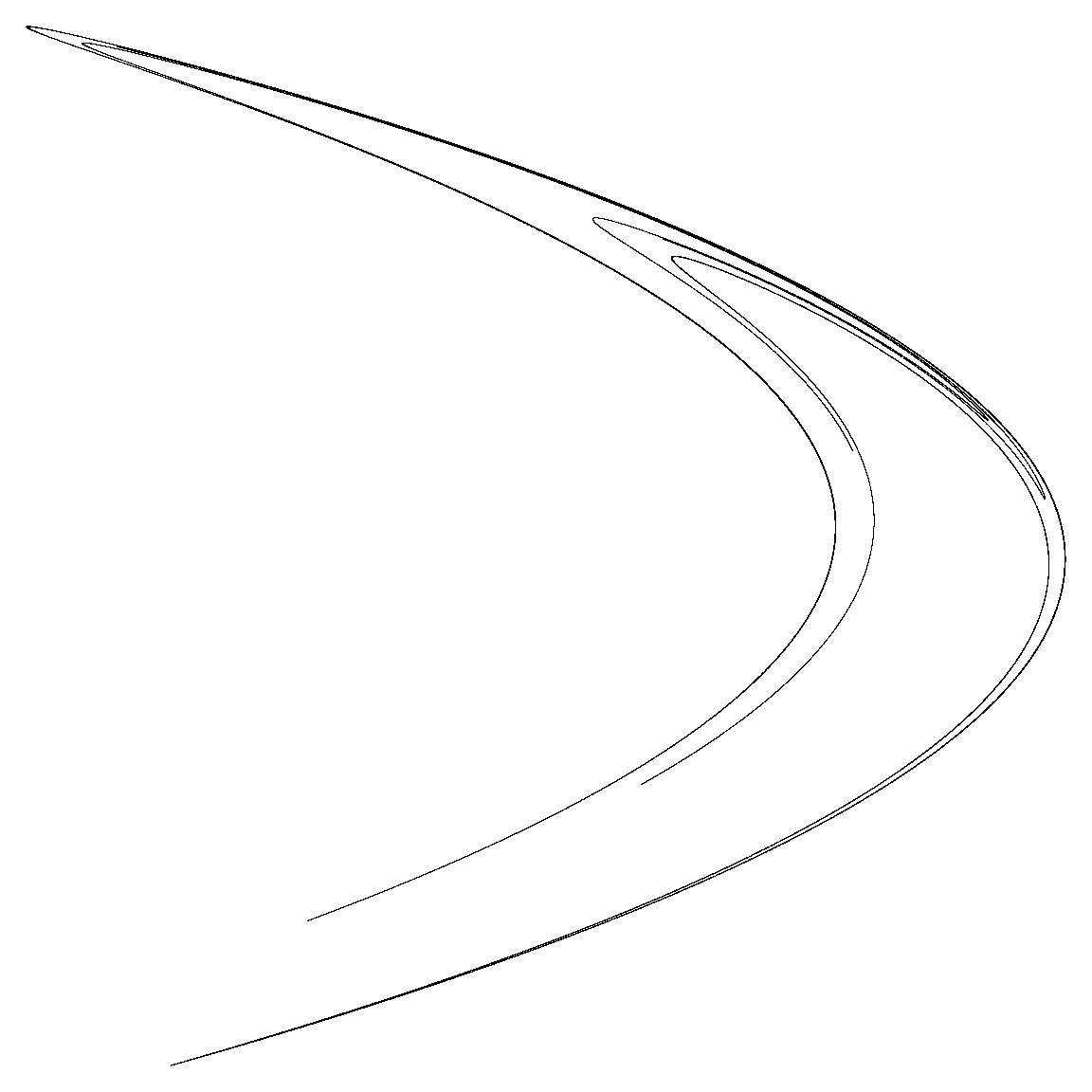}
\caption{Simulation of the H\'enon attractor for $a=1,4$ and $b=0,3$. Notice that the existing results
do not cover these `historical' values.}
\label{henon-fig}
\end{figure}

An important example of a dynamical system in the plane modeled by a Young tower with a return time decaying
exponentially is the Lozi map:
\[
T_{a,b}:
\begin{pmatrix}
x\\
y
\end{pmatrix}
\mapsto
\begin{pmatrix}
1-a|x|+y\\
bx
\end{pmatrix}
\]
which possesses an attractor depicted in Fig. \ref{A-lozi}. 
Lozi's map is much simpler to analyse than the famous
H\'enon map:
\[
T_{a,b}:
\begin{pmatrix}
x\\
y
\end{pmatrix}
\mapsto
\begin{pmatrix}
1-ax^2+y\\
bx
\end{pmatrix}
\]
For certain parameters, this map has an attactor displayed in Fig. \ref{henon-fig}.
For the so-called Benedicks-Carleson parameters\footnote{These parameters form a subset of $\mathds{R}^2$ with positive
Lebesgue measure \cite{BC}.}, it is possible to prove \cite{BY} that
the H\'enon attractor fits the general scheme of Young towers with
exponential tails. In particular, there is a unique SRB measure whose support is the
attractor.

\medskip

Important examples of maps, which are conservative, are billiard maps, like
planar Lorentz gases and Sinai's billiard. They can be also
modeled by Young towers. We refer to \cite{young1} but also
to \cite{CY} for a conceptual account avoiding technicalities.


\section{Limit theorems}
\label{sec:3}

In this section we review some limit theorems obtained for the class of systems previously described.

\subsection{Covariance and decay of correlations}


\begin{definition}[Correlations]\label{correl}
\begin{itemize}
\item[] For a dynamical system $(\Omega,T,\mu)$ and an observable $f:\Omega\to\mathds{R}$ in $L^2(\mu)$,
the \emph{autocovariance} of order $\ell\geq 0$ of the process $\{f\circ T^k; k\geq 0\}$ is defined as
\[
C_f(\ell):=\int f\cdot f\circ T^\ell \textup{d}\mu-\Big( \int f\textup{d}\mu\Big)^2.
\]
More generally, for a pair $f,g$ of observables in $L^2(\mu)$, the covariance  of order $\ell$ of
the processes $\{f\circ T^k; k\geq 0\}$ and $\{g\circ T^k; k\geq 0\}$ is defined as 
\[
C_{f,g}(\ell):=\int f\cdot g\circ T^\ell \textup{d}\mu-\int f\textup{d}\mu\int g\textup{d}\mu.
\]
\item[] In dynamical systems, it is customary to call the auto-covariance of order $\ell$
the ``\emph{correlation coefficient}'' of order $\ell$.
\end{itemize}
\end{definition}

The auto-covariance, or more generally, the covariance, is the basic indicator of a chaotic behavior:
for large values of $\ell$, the random variables $f$ and $f\circ T^\ell$ should be nearly independent,
{\em i.e.} the coefficient $C_f(\ell)$ should decay to $0$ as $\ell$ grows. Two factors affect the
rapidity of this decay: the strength of chaos in the underlying dynamical system $T:\Omega\circlearrowleft$
and the regularity of the observable $f$.

Recall that a dynamical system $(\Omega,\mathscr{B},T,\mu)$ is mixing if for any two measurable sets $A,B\subset \Omega$ one has
$\mu(A\cap T^{-n}B)\xrightarrow[n\to\infty]{}\mu(A)\mu(B)$. It is easy to prove that the system is mixing if and only
if correlations decay, {\em i.e.}, $C_{f,g}(\ell)\xrightarrow[n\to\infty]{} 0$ for every pair of $f,g\in L^2(\mu)$.

The speed or rate of the decay of correlations (also called the rate of mixing) is crucial in the statistical analysis of
chaotic systems.

\begin{svgraybox}
\begin{theorem}[Mixing and decay of correlations \cite{young1,young2,sarig,seb-sharp}]\label{mixing-towers}
\leavevmode\\
Let $T: \Omega\circlearrowleft$ be a dynamical system modeled by a Young tower and $\mu$ its SRB measure.
The system is mixing and the rate of decay of correlations for H\"older continuous observables is directly related to the behavior of $m^u\{R > n\}$ as $n\to\infty$.
\begin{itemize}
\item For example, if $m^u\{R > n\}= \mathcal{O}(e^{-an})$ for some $a>0$, then $(T,\mu)$ has exponential decay of correlations.
\item
If $m^u\{R > n\} = \mathcal{O}(1/n^\gamma)$ for some $\gamma>1$, then $(T,\mu)$ has polynomial decay of correlations.
More precisely, $C_f(\ell)=\mathcal{O}(1/\ell^{\gamma-1})$.
\end{itemize}
\end{theorem}
\end{svgraybox}
For the H\'enon map with Benedicks-Carleson parameters, correlations for H\"older continuous
observables decay exponentially fast. The intermittent map \eqref{MPmap} has polynomial
decay of correlations: $C_f(\ell)=\mathcal{O}(1/\ell^{\frac{1}{\alpha}-1})$.
Two-dimensional examples with an intermittent behavior come from billiards. Chernov
and Zhang studied in  \cite{CZ1,CZ2} several classes of billiards for which the decay of correlations
is $\mathcal{O}((\log \ell)^c/\ell^{1/\alpha -1})$ for some parameter $\alpha$ taking values
in $(0,1/2]$.

\subsection{Central limit theorem}\label{sec:clt}

We start by a definition.
\begin{definition}[Central limit theorem]
\begin{itemize}
\item[]Let $(\Omega,T,\mu)$ be a dynamical system  and $f:\Omega\to\mathds{R}$ an observable in $L^2(\mu)$.
We say that $f$ satisfies the central limit theorem (CLT for short) with respect to $(T,\mu)$ if there exists
$\sigma_f\geq 0$ such that
\begin{equation}\label{clt-ds}
\lim_{n\to\infty} \mu\left\{x : \frac{S_nf(x)-n\int f\textup{d}\mu}{\sqrt{n}}\leq t \right\}=
\frac{1}{\sqrt{2\pi}\sigma_f}\int_{-\infty}^t e^{-\frac{u^2}{2\sigma^2_f}} \textup{d}u,\quad \forall t\in\mathds{R}.
\end{equation}
In probabilistic notation, the previous convergence is written compactly as
\[
\frac{S_nf-n\int f\textup{d}\mu}{\sqrt{n}}\lawto \mathcal{N}_{0,\sigma^2_f},
\]
where $\mathcal{N}_{0,\sigma^2_f}$ stands for the Gaussian law with mean $0$ and
variance $\sigma^2_f$.\newline
When $\sigma_f=0$ the right-hand side has to be understood as the Heaviside function.
\end{itemize}
\end{definition}

In probabilistic terms, this definition asks for the \emph{convergence in law} of the ergodic average
`zoomed out' by the factor $\sqrt{n}$ to a random variable whose law is $\mathcal{N}_{0,\sigma^2_f}$.

By analogy with i.i.d. processes, one expects that $\sigma_f$ be the variance of the process $\{f\circ T^n\}$.
If it were an i.i.d. process, we would have
\[
\sigma_f^2=\textup{Var}\left(S_nf/\sqrt{n}\right)=
\int f^2\textup{d}\mu-\left(\int f\textup{d}\mu\right)^2\textup{d}\mu=C_f(0),
\]
where $\textup{Var}(X)=\mathds{E}\big[(X-\mathds{E}(X))^2\big]$ is the variance of $X$.
But because of the correlations between $f$ and $f\circ T^n$, this is not the case.
A natural candidate for the variance is
\[
\sigma_f^2=\lim_{n\to\infty}\frac{1}{n}\int \big(S_n f-n\int f\textup{d}\mu\big)^2 \textup{d}\mu,
\]
provided the limit exists.
Simple algebra, using the invariance of $\mu$ under $T$, gives
\[
\frac{1}{n}\int \big(S_n f-n\int f\textup{d}\mu\big)^2 \textup{d}\mu=C_f(0) + 2\sum_{\ell=1}^{n-1}\frac{n-\ell}{n}
C_f(\ell).
\]
It is simple to prove that if
\[
\sum_{j=1}^\infty |C_f(j)|<\infty
\]
then 
\[
\lim_{n\to\infty}\sum_{\ell=1}^{n-1}\frac{n-\ell}{n}\ C_f(\ell)=
\sum_{\ell=1}^\infty C_f(\ell),
\]
whence
\begin{equation}\label{variance}
\sigma_f^2=C_f(0)+2\sum_{\ell=1}^\infty C_f(\ell).
\end{equation}
We have the following theorem.
\begin{svgraybox}
\begin{theorem}[Central limit theorem, \cite{young1,young2}]\label{tclfortowers}
\leavevmode\\
Let $T:\Omega\circlearrowleft$ be a dynamical system modeled by a Young tower and $\mu$ its SRB measure.
Let $f:\Omega\to\mathds{R}$ be a H\"older continuous observable. If $\int R^2 \textup{d} m^u<\infty$ (which
implies $\sum_{\ell\geq 1}|C_f(\ell)|<\infty$), then $f$ satisfies the central limit theorem with respect to $(T,\mu)$.
\end{theorem}
\end{svgraybox}
For the class of systems discussed in this paper, it is well-known that typically
$\sigma_f^2>0$. Indeed, $\sigma_f^2=0$ only for H\"older observables lying in a closed subspace of infinite codimension.

For example, H\"older continuous observables satisfy the CLT for 
the H\'enon map with Benedicks-Carleson parameters. For the map
\eqref{MPmap}, the CLT holds if $\alpha<1/2$. We shall see what happens
when $\alpha\geq 1/2$ later on.\newline
There are examples of convergence to the Gaussian law but with a non-classical
renormalizing sequence $(\sqrt{n\log n})$, instead of $(\sqrt{n})$. This is the case
for Bunimovich's billiard (stadium) where correlations decay only as $1/n$ (where $n$ is
the number of collisions); see \cite{balint}.

\medskip

In essence, the central limit theorem tells us that {\em typically} ({\em i.e.} with very high probability),
\[
S_n f-n\int f\textup{d}\mu = \mathcal{O}(\sqrt{n}).
\]
In other words, the typical fluctuations of $S_n f/n$ around $\int f\textup{d}\mu$
are of order $1/\sqrt{n}$. But, in principle, $S_n f$ can take values as large
as $n$, {\em i.e.} $S_nf/n-\int f\textup{d}\mu$ can be of order one, but
with a small probability. Such fluctuations are naturally called `large deviations'.
This is the subject of the next section.

\subsection{Large deviations}\label{subsec:LD}

For a bounded i.i.d. process $\{X_n\}$, it is a classical result in probability, usually called Cram\'er's theorem \cite{DZ},
that
$
\mathds{P}\big\{ \big|n^{-1}\big(X_0+\cdots+X_{n-1}\big)-\mathds{E}[X_0]\big|>\delta\big\}
$
decays exponentially with $n$. Moreover,
\[
\lim_{n\to\infty}\frac{1}{n}
 \log \mathds{P}\left\{ \left|\frac{X_0+\cdots+X_{n-1}}{n}-\mathds{E}[X_0]\right|>\delta\right\}
= -\mathbf{I}(\delta).
\]
Typically, the function $\mathbf{I}$ (the so-called \emph{rate function}) is strictly convex and vanishes only
at $0$\footnote{The rate function must vanish at $0$ in view of Birkhoff's ergodic theorem.} (hence it is non-negative).
Since the process is bounded, its domain is a finite interval.
The rate function turns out to be the Legendre transform of the cumulant generating
function $\theta\mapsto \log \mathds{E}[\exp(\theta X_0)]$.\newline
One expects this exponential decay for the probability of deviation in `sufficiently chaotic' dynamical systems and for
a H\"older continuous observable $f$. For notational convenience, assume that $\int f\textup{d}\mu=0$.
The goal is to prove that there exists a rate function
$\mathbf{I}_f:\mathds{R}\to[0,+\infty]$ such that
\[
\lim_{\epsilon\to 0}\lim_{n\to\infty} \frac{1}{n}\log\mu\left\{x\in\Omega: \frac{1}{n} S_n f(x)\in [a-\epsilon,a+\epsilon]\right\}=-\mathbf{I}_f(a).
\]
In many situations, such a result is obtained by proving that the cumulant generating function
\[
\Psi_f(z)=\lim_{n\to\infty}\ \frac1n \log \int e^{z S_n f}\textup{d}\mu
\]
exists and is smooth enough for $z$ real in an interval containing the origin. Then the rate function
is the Legendre transform of $\Psi_f$. 
However, as we shall see, when chaos is not strong enough, one may indeed get subexponential decay rates for large deviations (and therefore there is no rate function).
  
For systems modeled by a Young tower with exponential tails, we have the following result. It
turns out that the logarithmic moment generating function $\Psi_f(z)$ can be studied for complex $z$.

\begin{svgraybox}
\begin{theorem}[Cumulant generating functions \cite{youngLD,melbourneLD}]\label{LMGF}
\leavevmode\\
Let $T: \Omega\circlearrowleft$ be a dynamical system modeled by a Young tower and $\mu$ its SRB measure.
Assume that  $m^u\{R > n\}= \mathcal{O}(e^{-an})$ for some $a>0$. Let $f:\Omega\to\mathds{R}$ be a
H\"older continuous observable such that $\int f \textup{d}\mu=0$. 
\begin{itemize}
\item Then there exist positive numbers $\eta=\eta(f)$ and $\xi=\xi(f)$ such that the
logarithmic moment generating function $\Psi_f$ exists and is analytic in the strip
\[
\{z\in\mathds{C}: |\textup{Re}(z)|<\eta,|\textup{Im}(z)|<\xi\}.
\]
\item In particular, $\Psi_f'(0)=\int f \textup{d}\mu$ and $\Psi"\!\!_f(0)=\sigma^2_f$,
which is the variance \eqref{variance} of the process $\{f\circ T^n\}$. Moreover, $\Psi_f(z)$
is strictly convex for real $z$ provided $\sigma_f^2>0$.
\end{itemize}
\end{theorem}
\end{svgraybox}
From this kind of result, one can deduce the following result by using
Gartner-Ellis theorem or the like (see \cite[section 4.5]{DZ} and \cite[pp. 102--103]{HH}).
Notice that it is enough for $\Psi_f$ to be differentiable to apply this theorem.

\begin{svgraybox}
\begin{theorem}[Exponential large deviations \cite{melbourneLD,youngLD}]\label{LDYexp}
\leavevmode\\
Under the same assumptions as in the previous theorem, let
$\mathbf{I}_f$ be the Legendre transform of
$\Psi_f$, {\em i.e.} $\mathbf{I}_f(t)=\sup_{z\in (-\eta,\eta)} \{tz-\Psi_f(z)\}$.
Then for any interval $[a,b]\subset [\Psi_f'(-\eta),\Psi_f'(\eta)]$,
\[
\lim_{n\to\infty} \frac1n \log\mu
\left\{x\in \Omega: \frac1n S_n f(x)\in [a,b]\right\}=-\inf_{t\in[a,b]} {\mathbf I}_f(t).
\]
\end{theorem}
\end{svgraybox}

\begin{remark}
Using a general theorem of Bryc \cite{bryc}, one can deduce the central limit theorem
from Theorem \ref{LMGF}. We stress that analyticity of $\Psi_f$ is necessary. In general,
if $\Psi_f$ is only $C^\infty$ (ensuring that $\Psi"\!\!_f(0)=\sigma^2_f$), it is false
than the central limit theorem follows from exponential large deviations. 
\end{remark}

We now turn to systems modeled by a Young tower with sub-exponential tails. In this case,
there is no rate function and one gets sub-exponential large deviation bounds.

\begin{svgraybox}
\begin{theorem}[Sub-exponential large deviations \cite{melbournePAMS}]\label{LDYpoly}
\leavevmode\\
Let $T: \Omega\circlearrowleft$ be a dynamical system modeled by a Young tower and $\mu$ its SRB measure.
Assume that  $m^u\{R > n\}= \mathcal{O}(1/n^\gamma)$ for some $\gamma>1$. Let $f:\Omega\to\mathds{R}$ be a H\"older continuous observable such that $\int f \textup{d}\mu=0$. Then, for any $\varepsilon>0$ 
\[
\mu\left\{x\in\Omega : \left| \frac{1}{n}S_n f(x)\right|>\varepsilon \right\}\leq
\frac{C_{f,\varepsilon}}{n^{\gamma-1}},\quad\textup{for any}\;n\in\mathds{N}.
\]
\end{theorem}
\end{svgraybox}

Notice that according to Theorem \ref{mixing-towers}, the decay is the same as that for
correlations. The dependence in $\varepsilon$ of the constant $C_{f,\varepsilon}$
is in $\varepsilon^{-2q}$ where $q>\max(1,\gamma-1)$.

Let us again use our favorite example, namely the Manneville-Pomeau map, to illustrate
the preceding result. In this case, one can also prove a lower bound for the
probability of large deviations. Indeed, for the map \eqref{MPmap}, the theorem applies
with $\gamma=\frac{1}{\alpha}$, where $\alpha\in (0,1)$. Recall that for $\alpha\in (0,1/2)$,
the central limit theorem holds (see Section \ref{sec:clt}), but it fails when $\alpha\in [1/2,1)$
(See Section \ref{subsec:nong} below).\newline
Moreover, it is proved in \cite{melbournePAMS} that 
there is a nonempty open set of H\"older observables $f$ for which $n^{-\frac{1}{\alpha}+1}$
is a lower bound for large deviations for $n$ sufficiently large. For these observables, we have
for any $\varepsilon>0$
\begin{equation}\label{borneinf}
\lim_{n\to\infty} \frac{\log \mu\left\{x\in [0,1]: \left|\frac{1}{n}S_n f(x)\right|>\varepsilon \right\}}{\log n}
=-\frac{1}{\alpha}+1.
\end{equation}

\subsection{Convergence to non-Gaussian laws}\label{subsec:nong}

The purpose of this section is to show what happens when the CLT fails but one
still has convergence in law, but with a renormalizing sequence different from
$(\sqrt{n})$.
For the reader's convenience, we recall the notion of domain of attraction
for an observable and a classical theorem about stable laws for i.i.d. processes.

A function $f$, defined on a probability space $(\Omega,\mathcal{B},m)$,
is said to \emph{belong to a domain of attraction} if it fulfills
one the following three conditions:

\begin{itemize}
\item[\textup{I}.] It belongs to $L^2(\Omega)$.
\item[\textup{II}.] One has $\int \mathds{1}_{\{|f|>x\}} \textup{d} m \sim
x^{-2}\ell(x)$, for some function $\ell$ such that $L(x):=2\int_1^x
\frac{\ell(u)}{u}\textup{d} u$ is of slow variation and unbounded.
\item[\textup{III}.] There exists $p\in (1,2)$
such that
\[
\int \mathds{1}_{\{f>x\}} \textup{d} m=(c_1+o(1))x^{-p}L(x)\quad\textup{and}
\quad \int \mathds{1}_{\{f<-x\}} \textup{d} m=(c_2+o(1))x^{-p}L(x),
\]
where $c_1,c_2$
are nonnegative real numbers such that $c_1+c_2>0$, and $L$ is of
slow variation.
\end{itemize}

Note that the three conditions are mutually exclusive.

The above definition of domain of attraction is motivated by the
following well-known, classical result in Probability (see 
{\em e.g.}~\cite{GK}):
\begin{svgraybox}
\begin{theorem}[Convergence to stable laws for i.i.d. processes]
\label{thm_limite_iid}
\leavevmode\\
Let $Z$ be a random variable belonging to a
domain of attraction. Let $Z_0,Z_1,\dots$ be a sequence of
independent, identically distributed, random variables with the same
law as $Z$. In all cases, we set $A_n=n\mathds{E}[Z]$ and
\begin{enumerate}
\item if condition I holds, we set $B_n=\sqrt{n}$ and $\mathcal{W}=
\mathcal{N}_{0, \mathds{E}[Z^2]-\mathds{E}[Z]^2}$;
\item if condition II holds, we let $B_n$ be a renormalizing sequence
with $nL(B_n)\sim B_n^2$, and $\mathcal{W}=\mathcal{N}_{0,1}$;
\item if condition III holds, we let $B_n$ be a renormalizing sequence
such that $nL(B_n)\sim B_n^p$. Define
$c=(c_1+c_2)\Gamma(1-p)\cos\left(\frac{p\pi}{2}\right)$ and
$\beta=\frac{c_1-c_2}{c_1+c_2}$.
\end{enumerate}
Let $\mathcal{W}=\mathcal{W}_{p,c,\beta}$ be the law with
characteristic function
\begin{equation}\label{cfstable}
\mathds{E}[e^{it\mathcal{W}}]=e^{-c|t|^p\left(1-i\beta \sgn(t)
\tan\left(\frac{p\pi}{2}\right) \right)},\quad (1<p\leq 2, c>0, |\beta|\leq 1).
\end{equation}
Then
\[ 
\frac{\sum_{i=0}^{n-1} Z_i -A_n}{B_n} \lawto \mathcal{W}.
\]
\end{theorem}
\end{svgraybox}

\bigskip

The case $p=2$ corresponds to the Gaussian law. For $p<2$, the corresponding distributions
are said to have `heavy tails' since $\mathds{P}\{Z>x\}=(c_1+o(1)) x^{-p}$ and
$\mathds{P}\{Z<-x\}=(c_2+o(1)) x^{-p}$.
The conditions put on the distribution of $Z$ are almost necessary
and sufficient to get a convergence in law of that type, we only
restricted the range of $p$'s, which could also be taken in the
interval $(0,1]$.

We illustrate the occurrence of non-Gaussian limit laws in the most important example, that is,
the Pomeau-Manneville map \eqref{MPmap}. 

\begin{svgraybox}
\begin{theorem}[Convergence to stable laws for the Manneville-Pomeau map \cite{seb-sharp}]\label{MPmapstable}
\leavevmode\\
Let $T_\alpha$ be the map of the interval \eqref{MPmap}, with $\alpha\in (0,1)$ and $\mu$ its
unique absolutely continuous, invariant, probability measure. Let $f:[0,1]\to\mathds{R}$ be a H\"older observable
and assume that $\int f \textup{d}\mu=0$.
\begin{itemize}
\item If $\alpha<1/2$ then the central limit theorem holds (this is a special instance of Theorem \ref{tclfortowers}).
\item If $\alpha>1/2$ then:
\begin{itemize}
\item if $f$ is Lipschitzian and $f(0)=0$, then the central limit theorem holds;
\item if $f(0)\neq 0$ then $\frac{1}{n^\alpha} S_n f$ converges in law to the stable law
$\mathcal{W}_{\frac{1}{\alpha},c,\sgn(f(0))}$ whose characteristic function is given by \eqref{cfstable}.
\end{itemize}
\end{itemize}
\end{theorem}
\end{svgraybox}

\medskip

When $\alpha=1/2$ and $f(0)\neq 0$, there is convergence to the Gaussian law but
with the unusual renormalizing sequence $(\sqrt{n\log n})$ (instead of $\sqrt{n}$).
See \cite{seb-CLT} for more details.

\subsection{Convergence in law made almost sure}\label{subsec:asclt}

The aim of this section is to show that whenever
we can prove a limit theorem in the classical sense for a dynamical
system, we can prove a suitable almost-sure version based on an
empirical measure with log-average.\newline
The prototype of such a theorem is the almost-sure central limit
theorem: if $X_n$ is an i.i.d.\ $L^2$ sequence with $\mathds{E}[X_i]=0$ and
$\mathds{E}[X_i^2]=1$, then, almost surely,
  \begin{equation}\label{ascltiid}
  \frac{1}{\log n} \sum_{k=1}^n \frac{1}{k} \gdelta_{\sum_{j=0}^{k-1}
  X_j/\sqrt{k}} \lawto \mathcal{N}_{0,1}
  \end{equation}
where ``$\lawto$'' means weak convergence of probability measures on
$\mathds{R}$. Here and henceforth, $\delta_x$ is the Dirac mass at $x$. This
result should be compared to the classical central limit theorem,
which can be stated as follows:
\[
\mathds{E}[ \mathds{1}_{\{\sum_{j=0}^{n-1} X_j/\sqrt{n}\
\leq t\}}]\xrightarrow[n\to\infty]{}
\frac{1}{\sqrt{2\pi}}\int_{-\infty}^ t e^{-u^2/2}\textup{d} u
\]
for any $t\in\mathds{R}$. To better compare these theorems, it is worth
noticing that \eqref{ascltiid} implies that \emph{almost surely}
  \begin{equation}
  \frac{1}{\log n} \sum_{k=1}^n \frac{1}{k}
  \mathds{1}_{\{\sum_{j=0}^{k-1} X_j/\sqrt{k}\ \leq t\}}
  \xrightarrow[n\to\infty]{}
  \frac{1}{\sqrt{2\pi}}\int_{-\infty}^ t e^{-u^2/2}\textup{d} u
  \end{equation}
for any $t\in\mathds{R}$. So, instead of taking the expected value, we take
a logarithmic average and obtain an almost-sure convergence.

In fact, whenever there is independence and a classical limit
theorem, the corresponding almost-sure limit theorem also holds
(under minor technical conditions), see \cite{berkes} and
references therein.

Let us put the following general definition:
\begin{definition}[Almost sure limit theorem towards a random variable]
\begin{itemize}
\item[]
Let $S_n$ be a sequence of random variables on a
probability space, and let $B_n$ be a renormalizing sequence.
\footnote{A renormalization function is a function $B: \mathds{R}_+^*\to
\mathds{R}_+^*$ of the form $B(x)=x^d L(x)$ where $d>0$ and $L$ is a
normalized slowly varying function. The corresponding
renormalizing sequence is $B_n:=B(n)$.} We say
that $S_n/B_n$ satisfies an almost sure limit theorem
towards a law $\mathcal{W}$ if, for almost all $\omega$,
\[
  \frac{1}{\log N} \sum_{k=1}^N \frac{1}{k} \gdelta_{ S_k(\omega)/B_k}
  \lawto \mathcal{W}.
\]
  \end{itemize}
\end{definition}
We now turn to the dynamical system context. The almost-sure central limit theorem, for instance,
takes the form
\[
\frac{1}{\log n} \sum_{k=1}^n \frac{1}{k}
\delta_{S_k f(x)/\sqrt{k}}\lawto\mathcal{N}_{0,\sigma_f^2},\;\textup{for}\; \mu-\textup{almost every}\;x,
\]
where, for notational simplicity, we assume that $\int f\textup{d}\mu=0$.

In the paper \cite{chazottes-gouezel}, we proved that ``whenever we can prove a limit
theorem in the classical sense for a dynamical system, we can prove
a suitable almost-sure version''.
More precisely, we investigated three methods that are used to prove limit theorems in
dynamical systems: spectral methods, martingale methods, and
induction arguments. We showed that whenever these methods apply,
the corresponding limit theorem admits a suitable almost-sure
version.\newline
For instance, one has the following result.
\begin{svgraybox}
\begin{theorem}[Convergence in law made almost sure \cite{chazottes-gouezel}]
\leavevmode\\
Let $T:\Omega\circlearrowleft$ be a dynamical system modeled by a Young tower
and let $\mu$ its SRB measure. Let $f:\Omega\to\mathds{R}$ be a H\"older continuous
observable such that $\int f \textup{d}\mu=0$. Then, if
\[
\quad\mu\left\{x\in\Omega: \frac{S_kf(x)}{B_k}\leq t\right\}\xrightarrow[k\to\infty]{} \mathcal{W}((-\infty,t]) 
\]
for every $t\in\mathds{R}$ at which $\mathcal{W}$ is continuous, for a certain law $\mathcal{W}$ and for a certain renormalizing sequence $(B_n)$, then
\[
  \frac{1}{\log n} \sum_{k=1}^n \frac{1}{k} \gdelta_{S_k f/ B_k}
  \lawto \mathcal{W} \quad\mu-\textup{almost-surely}.
\]
\end{theorem}
\end{svgraybox}

\medskip

Let us illustrate this theorem with a few examples. For any dynamical system modeled by
a Young tower with $L^2$ tails, one has
\[
\frac{1}{\log n} \sum_{k=1}^n \frac{1}{k}
\delta_{S_k f(x)/\sigma_f\sqrt{k}}\lawto\mathcal{N}_{0,\sigma_f^2}.
\]
For the Manneville-Pomeau map \eqref{MPmap}, this is true for $\alpha\in(0,1/2)$. When
$\alpha>1/2$, this is still the case provided that $f(0)=0$ and $f$ is Lipschitz. If $f(0)\neq 0$, then 
\[
\frac{1}{\log n} \sum_{k=1}^n \frac{1}{k}
\delta_{S_k f(x)/k^\alpha}\lawto\mathcal{W}_{\frac{1}{\alpha},c,\sgn(f(0))}
\]
(see Theorem \ref{MPmapstable}).

\section{Concentration inequalities and applications}
\label{sec:4}

\subsection{Introduction}

We start by the simplest occurrence of the concentration of measure phenomenon \cite{ledoux}.
Consider an independent sequence of Bernoulli random variables $(\eta_i)_{0\leq i\leq n-1}$
({\em i.e.} $\mathds{P}(\eta_i=-1)=\mathds{P}(\eta_i=1)=1/2$, whence $\mathds{E}[\eta_i]=0$).
Then one has the following classical inequality (Chernov's bound):
\begin{equation}\label{chernov}
\mathds{P}\Big( \Big|\small{\sum_{i=0}^{n-1} \eta_i}\Big|\geq t\Big)\leq
2\exp\Big(-\frac{t^2}{2n}\Big),\quad \forall t\geq 0.
\end{equation}
This exponential inequality reflects the most important theorem of probability, imprecisely stated as follows: 
``In a long sequence of tossing a fair coin, it is likely that heads will come up nearly half of the time.''
Indeed, if we let $B_n$ be the number of $1$'s in the sequence $(\eta_i)_{0\leq i\leq n-1}$, then
$\sum_{i=0}^{n-1} \eta_i=2B_n-n$, and so \eqref{chernov} is equivalent to
\[
\mathds{P}\left( \Big|B_n-\frac{n}{2}\Big|\geq t\right)\leq 2\exp\Big(\frac{-2t^2}{n}\Big),\quad \forall t\geq 0.
\]
This is of course a much stronger statement than the Strong Law of Large Numbers.

The perspective of concentration inequalities is to look at the random variable
$Z_n=\sum_{i=0}^{n-1} \eta_i$ as a function of the individual variables $\eta_i$. 
Inequality \eqref{chernov}, when $Z_n$ is normalized by $n$ (since it can take values as large as $n$)
can be phrased pretty offensively by saying that
\begin{center}
$\frac{Z_n}{n}$ is essentially constant ($=0$).
\end{center}
The scope of concentration inequalities is to understand to what extent a general function $K$ of
$n$ random variables $X_0,\ldots,X_{n-1}$, and not just the sum of them, 
concentrates around its expectation like a sum of Bernoulli random variables. Of course,
the smoothness of $K$ has to play a role, as well as the dependence between the $X_i$'s.

Stated informally as a principle, the measure of concentration phenomenon is the following:
\begin{itemize}
\item[]``A random variable that smoothly depends on the influence of many weakly dependent
random variables is, on the appropriate scale, very close to a constant.''
\end{itemize}
This statement is of course quantified by statements like \eqref{chernov} or weaker ones, as we shall see.

In the context of dynamical systems, there are many examples of random variables
$K(X_0,\ldots,X_{n-1})$ which appear naturally but are defined in an indirect or complicated way.
Concentration inequalities, when available, allow to obtain, in a systematic way, {\em a priori}
bounds on the fluctuations of $K(X_0,\ldots,X_{n-1})$ around its expectation by using a simple information on $K$, namely
its Lipschitz constants.

\subsection{Concentration inequalities: abstract definitions}

We formulate some abstract definitions.

Let $\Omega$ be a metric space. A real-valued function $K$ on $\Omega^n$ is separately Lipschitz if, for any $i$, there
exists a constant $\textup{Lip}_i(K)$ such that
 \begin{multline*}
\big| K(x_0,\ldots,x_{i-1},x_i,x_{i+1},\ldots,x_{n-1})-K(x_0,\ldots,x_{i-1},x'_i,x_{i+1},\ldots,x_{n-1}) \big| \\
 \leq \textup{Lip}_i(K) d(x_i,x'_i)
\end{multline*}
for all points $x_0,\ldots,x_{n-1},x'_i$ in $\Omega$.\newline
Consider a stationary process $\{X_0,X_1,\ldots\}$ taking values in $\Omega$.

\begin{definition}[Exponential concentration inequality]
\begin{itemize}
\item[]
We say that the process  $\{X_0,X_1,\ldots\}$
satisfies an exponential concentration inequality if there exists a constant $C>0$ such that, for any separately
Lipschitz function $K(x_0,\ldots,x_{n-1})$, one has
\begin{equation}\label{exp-ineq}
\mathds{E}\left[ e^{K(X_0,\ldots,X_{n-1})-\mathds{E}[K(X_0,\ldots,X_{n-1})]}\right]
\leq e^{C \sum_{\ell=0}^{n-1}\textup{Lip}_\ell(K)^2}.
\end{equation}
\end{itemize}
\end{definition}

\medskip

In some cases, it is not reasonable to hope for such a strong inequality. This leads to the following
definition.
\begin{definition}[Polynomial concentration inequality]\label{moment-ineq}
\begin{itemize}
\item[] We say that the process  $\{X_0,X_1,\ldots\}$
satisfies a polynomial concentration inequality with moment $p\geq 2$  if there exists a constant $C>0$ such that, for any separately Lipschitz function $K(x_0,\ldots,x_{n-1})$, one has
\begin{equation}\label{poly-ineq}
\mathds{E}\left[\left|K(X_0,\ldots,X_{n-1})-\mathds{E}[K(X_0,\ldots,X_{n-1})]\right|^p\right]
\leq C \left(\sum_{\ell=0}^{n-1}\textup{Lip}_\ell(K)^2\right)^{p/2}.
\end{equation}
\end{itemize}
\end{definition}

\medskip


An important special case of \eqref{poly-ineq} is for $p=2$, which gives an inequality
for the variance of $K(X_0,\ldots,X_{n-1})$:
\begin{equation}\label{devroye}
\textup{Var}\big(K(X_0,\ldots,X_{n-1})\big)\leq C\ \sum_{\ell=0}^{n-1}\textup{Lip}_\ell(K)^2.
\end{equation}

\medskip

After these definitions, a few comments are in order.
\begin{itemize}
\item The crucial point in  \eqref{exp-ineq} and  \eqref{poly-ineq} is that the constant $C$
does depends neither on $K$ nor on $n$. It solely depends on the process.
\item These inequalities are not asymptotic, they hold true for any $n$.
\item  Obviously \eqref{exp-ineq} is a much stronger inequality than \eqref{poly-ineq}.
For instance, one can get \eqref{devroye} from \eqref{exp-ineq} as follows:
Multiply $K$ by $\lambda\neq 0$, substract $1$ from both sides, divide by $\lambda^2$; conclude
by using Taylor expansion and by letting $\lambda$ go to $0$.
\item An important consequence of the previous inequalities is a control on the
deviation probabilities of $K(X_0,\ldots,X_{n-1})$ from its expectation:

If a stationary process $\{X_n\}$ satisfies the exponential concentration inequality \eqref{exp-ineq}
then, for any $t>0$, one has
\begin{equation}\label{exp-dev}
\mathds{P}\left\{\left| K(X_0,\ldots,X_{n-1})-\mathds{E}[K(X_0,\ldots,X_{n-1})]\right|>t \right\}\leq 
2\ e^{-\frac{t^2}{4C\sum_{\ell=0}^{n-1}\textup{Lip}_\ell(K)^2}}.
\end{equation}
If the process satisifies the polynomial concentration inequality \eqref{poly-ineq}, one gets that for any $t>0$
\begin{equation}\label{poly-dev}
\mathds{P}\left\{\left| K(X_0,\ldots,X_{n-1})-\mathds{E}[K(X_0,\ldots,X_{n-1})]\right|>t \right\}\leq 
C t^{-q}\left(\sum_{\ell=0}^{n-1}\textup{Lip}_\ell(K)^2\right)^{q/2}.
\end{equation}
To prove \eqref{exp-dev}, we use Markov's inequality and \eqref{exp-ineq}: for any $t,\lambda>0$
\begin{align*}
& \mathds{P}\left\{K(X_0,\ldots,X_{n-1})-\mathds{E}[K(X_0,\ldots,X_{n-1})]>t \right\} \\
& \qquad =
\mathds{P}\left\{\exp\Big(\lambda \big(K(X_0,\ldots,X_{n-1})-\mathds{E}[K(X_0,\ldots,X_{n-1})]\big)\Big)>\exp(\lambda t) \right\} \\
&\qquad \leq e^{-\lambda t}
\mathds{E}\left[ e^{\lambda \big(K(X_0,\ldots,X_{n-1})-\mathds{E}[K(X_0,\ldots,X_{n-1})]\big)}\right]\\
& \qquad\leq e^{-\lambda t} e^{C \lambda^2\sum_{\ell=0}^{n-1}\textup{Lip}_\ell(K)^2}.
\end{align*}
This upper bound is minimized when $\lambda=t/(2C\sum_{\ell=0}^{n-1}\textup{Lip}_\ell(K)^2)$, whence
\[
\mathds{P}\left\{K(X_0,\ldots,X_{n-1})-\mathds{E}[K(X_0,\ldots,X_{n-1})]>t \right\}\leq 
e^{-\frac{t^2}{4C\sum_{\ell=0}^{n-1}\textup{Lip}_\ell(K)^2}}.
\]
The previous procedure is usually called the `Chernoff bounding trick'.
Of course, we can apply this inequality to $-K$ and deduce at once \eqref{exp-dev}.\newline
Inequality \eqref{poly-dev} follows immediately from Markov's inequality.\hfill $\square$
\end{itemize}

\subsection{Concentration inequalities for dynamical systems}

We now present concentration inequalities in the setting of non-uniformly hyperbolic
dynamical systems. In a forthcoming paper with S. Gou\"ezel \cite{jrseb} we prove the following theorems.
Let us notice that we take separately Lipschitz observables for the sake
of simplicity. All results are valid in the H\"older case (see \cite[Section 7.1]{jrseb}).

\subsubsection{Main results}

\begin{svgraybox}
\begin{theorem}[Exponential concentration inequality \cite{jrseb}]\label{exp-jrseb}
\leavevmode\\
Let $(\Omega,T,\mu)$ be a dynamical system modeled by a Young tower with exponential tails.
Then it satisfies an exponential concentration inequality: there exists a constant $C>0$ such that, for any
$n\in\mathds{N}$, for any separately Lipschitz function $K(x_0,\ldots,x_{n-1})$, 
\begin{equation}\label{exp-ineq-ds}
\int e^{K(x,Tx,\ldots,T^{n-1}x)-\int K(y,Ty,\ldots,T^{n-1}y)\textup{d}\mu(y)}
\textup{d}\mu(x) \leq
e^{C \sum_{\ell=0}^{n-1}\textup{Lip}_\ell(K)^2}
\end{equation}
\end{theorem}
\end{svgraybox}

As a consequence of the Chernoff bounding trick (see the previous section),
we get, for any $t>0$ and for any $n\in\mathds{N}$,
 \begin{multline}\label{devexp-ds}
 \mu\left\{ x \in \Omega: K(x,Tx,\dotsc, T^{n-1}x) -\int K(y,\dotsc, T^{n-1}y)\textup{d} \mu(y)>t\right\}\\
 \leq e^{-\frac{t^2}{4C\sum_{j=0}^{n-1} \textup{Lip}_j(K)^2}}.
 \end{multline}
The same bound holds for lower deviations by applying~\eqref{devexp-ds}
to $-K$.

There are well-known dynamical systems $(X,T)$ which can be modeled
by a Young tower with exponential tails
\cite{young1}. Examples of invertible dynamical systems
fitting this framework are for instance Axiom A attractors, H\'enon's
attractor for Benedicks-Carleson parameters \cite{BY},
piecewise hyperbolic maps like the Lozi attractor, some billiards
with convex scatterers, etc. A non-invertible example is the
quadratic family for Benedicks-Carleson parameters.

\begin{svgraybox}
\begin{theorem}[Polynomial concentration inequality \cite{jrseb}]\label{poly-jrseb}
\leavevmode\\
Let $(\Omega,T,\mu)$ be a dynamical system modeled by a Young tower. Assume that, for some
$q\geq 2$, $\int R^q \textup{d}m^u<\infty$. Then it satisfies a polynomial concentration inequality with
moment $2q-2$, {\em i.e.}, there exists a constant $C>0$ such that, for any
$n\in\mathds{N}$, for any separately Lipschitz function $K(x_0,\ldots,x_{n-1})$, 
\begin{multline}\label{poly-ineq-ds}
\int \left|K(x,Tx,\ldots,T^{n-1}x)- \int K(y,Ty,\ldots,T^{n-1}y)\textup{d}\mu(y)\right|^{2q-2}
\textup{d}\mu(x)
\leq\\ C \left(\sum_{\ell=0}^{n-1}\textup{Lip}_\ell(K)^2\right)^{q-1}.
\end{multline}
\end{theorem}
\end{svgraybox}

As a direct application of Markov's inequality, we get from that, for any $t>0$ and for any $n\in\mathds{N}$,
\begin{multline}\label{poly-dev-ds}
\mu\left\{ x \in \Omega: \big|K(x,Tx,\dotsc, T^{n-1}x) -\int K(y,\dotsc, T^{n-1}y)\textup{d} \mu(y)\big|>t\right\}
 \leq \\
C\  \frac{\left(\sum_{\ell=0}^{n-1}\textup{Lip}_\ell(K)^2\right)^{q-1}}{t^{2q-2}}
\end{multline}

For the Manneville-Pomeau map, we know that the exponential concentration inequality cannot be true. Indeed, 
\eqref{borneinf} is clearly an obstruction. Applying Theorem \ref{poly-jrseb}, we get a concentration inequality 
with moment $Q$ for any $Q<\frac{2}{\alpha}-2$ when $\alpha\in (0,1/2)$. Applying \eqref{poly-dev} yields
a deviation bound in $n^{-\frac{1}{\alpha}+1+\delta}$, for any $\delta>0$. This is very close to the 
upper bound in $n^{-\frac{1}{\alpha}+1}$ guaranteed by Theorem \ref{LDYpoly}. In fact, one can get
an optimal deviation inequality and get the latter bound, but we need the notion of a weak polynomial concentration inequality
that we do not want to detail here, see \cite{jrseb}. 

\subsubsection{About the literature}

The first paper in which a concentration inequalities was proved for dynamical systems
is \cite{cms}: an exponential concentration inequality is established for piecewise
uniformly expanding maps of the interval. 
For dynamical systems $(X,T)$ modeled by a Young tower with exponential tails,
a polynomial concentration inequality with moment $2$ (variance) was proved
in  \cite{CCS1}. 
Regarding systems with subexponential decay of correlations, the first result
was obtained in \cite{CCRV} for the Manneville-Pomeau map \eqref{MPmap}:
a polynomial concentration inequality with moment $2$ was proved for
$\alpha\leq 4-\sqrt{15}$.
The above theorems, proved in \cite{jrseb}, improve all these results in several ways.



\subsection{A sample of applications of concentration inequalities}\label{subsec:applications}

We present some applications of concentration inequalities to show them
in action. Some more, as well as all proofs, can be found in \cite{CCS2,chazottes-gouezel,CM,cms}.

\subsubsection{Warming-up with ergodic sums}

Let us apply the exponential inequality to the basic example is
$K_0(x_0,\ldots,x_{n-1})=f(x_0)+\cdots +f(x_{n-1})$ where  $f$ is a Lipschitz observable.
We obviously have $\textup{Lip}_i(K_0)=\textup{Lip}(f)$ for any $i=0,\ldots,n-1$.
When evaluated along an orbit segment $x,\ldots,T^{n-1}x$, we of course get the ergodic sum $S_nf(x)$.
Assuming that \eqref{devexp-ds} holds one gets
\[
\mu\left\{x\in\Omega :  \left|\frac{1}{n}S_n f(x)-\int f\textup{d}\mu\right|>t\right\}\leq 2 e^{-\frac{ nt^2}{4C\textup{Lip}(f)^2}},\;\forall t>0.
\]
Compared with large deviations (see Subsection \ref{subsec:LD}), we observe that this is the right order in $n$.
The large deviation result provides a much more accurate description of this deviation probability as $n\to\infty$.
But the previous inequality shows how small this deviation probability is already for finite $n$'s.

\subsubsection{Correlations}

Let $(\Omega,T,\mu)$ be an ergodic dynamical system and  $f:\Omega\to\mathds{R}$ be
a Lipschitz observable such that $\int f\textup{d}\mu=0$.
An obvious estimator of the correlation coefficient  $C_f(k)$ (cf. Def. \ref{correl}) is
\[
\widehat{C}_f(n,k,x)=\frac{1}{n} \sum_{j=0}^{n-1} f(T^j x)f(T^{j+k} x).
\]
Indeed, an immediate consequence of Birkhoff's ergodic theorem is that
\[
\widehat{C}_f(n,k,x)\xrightarrow[n\to\infty]{} C_f(k),\quad \mu-\textup{a.s.}
\]
Observe that $\int \widehat{C}_f(n,k,x) \textup{d}\mu= C_f(k)$ by the invariance of the measure.

We have the following result.

\begin{svgraybox}
\begin{theorem}[Correlation coefficients]
\leavevmode\\
Let $T:\Omega\circlearrowleft$ be a dynamical system modeled by a Young tower and $\mu$ its
SRB measure.
Let $f:\Omega\to\mathds{R}$ be a Lipschitz observable such that $\int f\textup{d}\mu=0$. 
\begin{itemize}
\item If the tower has exponential tails, there exists
$D>0$ such that for any $t>0$ and any $k,n\in\mathds{N}$
\[
\mu\left\{ x\in\Omega : \left|\widehat{C}_f(n,k,x)-C_f(k) \right|>t\right\}\leq 2 e^{-D \frac{n^2 t^2}{n+k}}.
\]
\item If, for some $q\geq 2$, $\int R^q \textup{d}m^u<\infty$, then there exists $G>0$ such that for any $t>0$
and any $k,n\in\mathds{N}$
\[
\mu\left\{ x\in\Omega : \left|\widehat{C}_f(n,k,x)-C_f(k) \right|>t\right\}\leq G \left(\frac{n+k}{n^2} \right)^{q-1} \frac{1}{t^{2q-2}}.
\]
\end{itemize}
\end{theorem}
\end{svgraybox}

The proof is easy. One considers the function
\[
K(x_0,\ldots,x_{n+k-1})=\frac{1}{n} \sum_{j=0}^{n-1} f(x_j)f(x_{j+k})
\]
of $n+k$ variables. It is obvious that $\textup{Lip}_i(K)\leq \|f\|_\infty \textup{Lip}(f)/n$. Applying
\eqref{exp-dev} and \eqref{poly-dev} yields the desired inequality.

\subsubsection{Empirical measure}

Let $(\Omega,T,\mu)$ be an ergodic dynamical system.
Birkhoff's ergodic theorem (see Subsection \ref{subsec:birk}) implies that the empirical measure
$\mathcal{E}_n(x)=(1/n) \sum_{j=0}^{n-1}\gdelta_{T^j x}$ converges
vaguely to $\mu$. We want to obtain a `speed' for this convergence, so we
need to define a distance. We use the Kantorovich distance $\textup{dist}_{\scriptscriptstyle{K}}$.
For two probability measures $\mu_1$ and $\mu_2$ on $\Omega$, it is defined as
\[
\textup{dist}_{\scriptscriptstyle{K}}(\mu_1,\mu_2)=\sup\left\{ \int g \textup{d}\mu_1-\int g \textup{d}\mu_2 : g:
\Omega\to\mathds{R}\;\textup{is}\;1-\textup{Lipschitz}\right\}.
\]
This distance is compatible with the vague topology.\newline
We are led to consider the observable
\[
K(x,Tx,\ldots,T^{n-1}x)=\textup{dist}_{\scriptscriptstyle{K}}\big(\mathcal{E}_n(x),\mu\big).
\]
\begin{svgraybox}
\begin{theorem}[Empirical measure]\label{emp}
\leavevmode\\ 
Let $T:\Omega\circlearrowleft$ be a dynamical system modeled by a Young tower with exponential tails
and $\mu$ its SRB measure. Then, for any $t>0$ and for any $n\in\mathds{N}$
\[
\mu\left\{x\in\Omega :
\Big|\textup{dist}_{\scriptscriptstyle{K}}(\mathcal{E}_n(x),\mu)-
\int \textup{dist}_{\scriptscriptstyle{K}}(\mathcal{E}_n(y),\mu)\textup{d}\mu(y)\Big| >\frac{t}{\sqrt{n}}\right\}
\leq 2 e^{-t^2/4C}.
\]
\end{theorem}
\end{svgraybox}

This theorem follows at once from \eqref{exp-dev} and the fact that the function $K$ defined above
has all its Lipschitz constants bounded by $1/n$. A natural step further is to try to get an upper
bound for $\int \textup{dist}_{\scriptscriptstyle{K}}(\mathcal{E}_n(\cdot),\mu)\textup{d}\mu$.
There is no general good bound in general; one has first to restrict to one-dimensional
systems (because there is a special representation for the Kantorovich distance in terms
of the distribution functions). Second, the regularity of the observables for which there is exponential decay of
correlations is crucial. We mention only one result for the quadratic map $T_a(x)=1-ax^2$ acting on
$\Omega=[-1,1]$, where $a\in [0,2]$. For Benedicks-Carleson parameters, we mentioned above
that this system can be modeled by a Young tower with exponential tails. In fact there is an 
exponential decay of correlations for more general observables than the H\"older ones, namely
for observables with bounded variation \cite{young4}. This allows to prove that
\[
\int \textup{dist}_{\scriptscriptstyle{K}}(\mathcal{E}_n(\cdot),\mu)\textup{d}\mu\leq \frac{B}{\sqrt{n}}
\]
for some $B>0$. Hence we deduce the following result from \eqref{emp}.

\begin{svgraybox}
\begin{theorem}
\leavevmode\\
Consider the map $T_a(x)=1-ax^2$ acting on $\Omega=[-1,1]$ for $a$ in the Benedicks-Carleson set
of parameters. Then there exist $D,t_0>0$ such that for any $t\geq t_0$ and for any $n\in\mathds{N}$
\[
\mu\left\{x\in\Omega :
\textup{dist}_{\scriptscriptstyle{K}}(\mathcal{E}_n(x),\mu)>\frac{t}{\sqrt{n}}\right\}
\leq 2 e^{-Dt^2}.
\]
\end{theorem}
\end{svgraybox}

\bigskip

A natural question is to estimate the density of the absolutely continuous invariant measure
of a one-dimensional dynamical system. A classical estimator is the so-called kernel
density estimator. We refer to \cite{CCS2,jrseb} for details and results.

\subsubsection{Tracing orbits}

We use concentration inequalities to quantify the tracing properties of some subsets of orbits. The basic problem can be formulated as follows. Let $A$ be a set of initial conditions and $x$ an initial condition not in $A$:
How well can one approximate the orbit of $x$ by an orbit from an initial condition of $A$ ?
One can measure the `average quality of tracing' by defining
\[
\mathcal{S}_A(x,n)= \frac{1}{n} \inf_{y\in A} \sum_{j=0}^{n-1} d(T^j x,T^j y)
\]
where $d$ is the distance on $\Omega$. Assume that $\textup{diam}(\Omega)=1$. We have the following result.
\begin{svgraybox}
\begin{theorem}
\leavevmode\\
Let $T:\Omega\circlearrowleft$ be a dynamical system modeled by a Young tower with exponential tails
and $\mu$ its SRB measure.
There exists a constant $c>0$ such that for any subset $A\subset X$ with strictly positive
$\mu$-measure, for any $n\in\mathds{N}$ and for any $t>0$ 
\[
\mu\left\{
x \in\Omega : \mathcal{S}_A(x,n)> c\frac{\sqrt{\log n}}{\mu(A)\sqrt{n}} + \frac{t}{\sqrt{n}}
\right\}
\leq e^{-t^2/4C}
\]
(where $C>0$ is the constant appearing in Theorem \ref{exp-jrseb}).
\end{theorem}
\end{svgraybox}

\begin{proof}
The function of $n$ variables
\[
K(x_0,\ldots,x_{n-1})=\frac{1}{n} \inf_{y\in A}\sum_{j=0}^{n-1} d\big(x_j,T^j y\big).
\]
is separately Lipschitz and it is easy to check that $\textup{Lip}_i(K)\leq 1/n$ for any $i=0,\ldots,n-1$. 
We use \eqref{exp-dev} to get at once
\begin{equation}\label{shadow}
\mu\left\{ x  : \mathcal{S}_A(x,n)>\int \mathcal{S}_A(y,n)\textup{d}\mu(y)+\frac{t}{\sqrt{n}}\right\}\leq
e^{-t^2/4C}.
\end{equation}
We now estimate $\int \mathcal{S}_A(y,n)\textup{d}\mu(y)$ from above. Fix $s>0$
and define the set
\[
B_s=\left\{ x : \mathcal{S}_A(x,n)> \int \mathcal{S}_A(y,n)\textup{d}\mu(y) + \frac{s}{\sqrt{n}}\right\}\cdot
\]
We have the identity
\[
\int \mathcal{S}_A(y,n)\textup{d}\mu(y) = 
\int_A \mathcal{S}_A(y,n)\textup{d}\mu(y)+\int_{A^c\cap B_s^c} \mathcal{S}_A(y,n)\textup{d}\mu(y)+
\int_{B_s} \mathcal{S}_A(y,n)\textup{d}\mu(y).
\]
The first integral is equal to $0$ by the very definition of $ \mathcal{S}_A$. The second one is bounded by
\[
\left(\int \mathcal{S}_A(y,n)\textup{d}\mu(y)+\frac{s}{\sqrt{n}}\right) \mu(A^c).
\]
And the third one is bounded by $\mu(B_s)$ because $\mathcal{S}_A(y,n)\leq 1$. By \eqref{shadow}
one has
\[
\mu(B_s)\leq e^{-s^2/4C}. 
\]
Hence
\[
\int \mathcal{S}_A(y,n)\textup{d}\mu(y)\leq
\left(\int \mathcal{S}_A(y,n)\textup{d}\mu(y)+\frac{s}{\sqrt{n}}\right) \mu(A^c) + e^{-s^2/4C},
\]
{\em i.e.}
\[
\int \mathcal{S}_A(y,n)\textup{d}\mu(y)\leq \mu(A)^{-1}\left(\frac{s}{\sqrt{n}} + e^{-s^2/4C}\right).
\]
To finish the proof, it remains to optimize over $s>0$.\hfill $\square$
\end{proof}

\medskip

For a system modeled by a Young tower with polynomial tails, one can obtain a weaker bound, see
\cite{jrseb}.

\subsubsection{Integrated periodogram}

Let $(\Omega,T,\mu)$ be an ergodic dynamical system and $f:\Omega\to\mathds{R}$ be a Lipschitz
observable with $\int f\textup{d}\mu=0$. Define the empirical integrated periodogram of
the process $\{f\circ T^k\}$ by 
\[
\EuScript{J}_n(x,\omega)=\int_0^\omega
\frac{1}{n} \Big| \sum_{j=0}^{n-1} e^{-\mathbf{i}js} f(T^j x)\Big|^2 \textup{d}s, \quad \omega\in[0,2\pi].
\]
Let
\[
\EuScript{J}(\omega)=C_f(0)\omega + 2\sum_{k=1}^\infty \frac{\sin(\omega k)}{k}\ C_f(k),
\]
that is, the cosine Fourier transform of the sequence of correlation coeficients. (Recall
that $C_f(k)=\int f\cdot f\circ T^k \textup{d}\mu$.)
One can prove the following theorem.

\begin{svgraybox}
\begin{theorem}
\leavevmode\\
Let $T:\Omega\circlearrowleft$ be a dynamical system modeled by a Young tower  with exponential tails
and $\mu$ its SRB measure. Let $f:\Omega\to\mathds{R}$ be a Lipschitz function such that $\int f\textup{d}\mu=0$.
There exist some positive constants $c_1,c_2$ such that
for any $n\in\mathds{N}$ and for any $t>0$
\[
\mu\left\{x \in \Omega: \sup_{\omega\in[0,2\pi]}\big|\EuScript{J}_n(x,\omega)-\EuScript{J}(\omega)\big|> t +
\frac{c_1(1+\log n)^{3/2}}{\sqrt{n}}\right\}
\leq e^{-c_2 n t^2/(1+\log n)^2}.
\]
\end{theorem}
\end{svgraybox}

The proof can be found in \cite{jrseb}.

\subsubsection{Almost-sure central limit theorem}

We come back to the almost-sure central limit theorem (cf. Subsection \ref{subsec:asclt}). Let $f$ be
a Lipschitz observable such that $\int f\textup{d}\mu=0$. For convenience, let
\[
\mathcal{A}_n=\frac{1}{\log n} \sum_{k=1}^n \frac{1}{k} \gdelta_{S_kf/\sqrt{k}}.
\]
This is a random measure on $\mathds{R}$. Given $x\in\Omega$, $\mathcal{A}_n(x)$ is a measure.
To measure its closeness to the Gaussian law
$\mathcal{N}_{0,\sigma_f^2}$, we use the Kantorovich distance $\textup{dist}_{\scriptscriptstyle{K}}$.
For two probability measures $\mu_1$ and $\mu_2$ on $\mathds{R}$, it is defined as
\[
\textup{dist}_{\scriptscriptstyle{K}}(\mu_1,\mu_2)=\sup\left\{ \int g \textup{d}\mu_1-\int g \textup{d}\mu_2 : g:
\mathds{R}\to\mathds{R}\;\textup{is}\;1-\textup{Lipschitz}\right\}.
\]
Convergence in this distance entails both weak convergence and convergence of the first moment.

\begin{svgraybox}
\begin{theorem}[Almost-sure central limit theorem]
\leavevmode\\
Let $T:\Omega\circlearrowleft$ be a dynamical system modeled by a Young
tower such that
\[
\int R^2 \textup{d}m^u<\infty
\] 
and $\mu$ its SRB measure. Let $f:\Omega\to\mathds{R}$ be a Lipschitz observable
with $\int f\textup{d}\mu=0$. Assume that $\sigma_f^2>0$. Then
\[
\textup{dist}_{\scriptscriptstyle{K}}\big(\mathcal{A}_n(x),\mathcal{N}_{0,\sigma_f^2}\big)\to 0\quad\textup{for}\;\mu-\textup{a.e.}\;x\in\Omega.
\]
\end{theorem}
\end{svgraybox}

\medskip

This is slightly stronger than the usual almost-sure central limit theorem. In fact, a more general statement is true:
if a process $\{X_k\}$ satisfies the central limit theorem and \eqref{devroye}, then the previous theorem is true.
This is the way it is proved in \cite{CCS2}.

\section{Open questions}
\label{sec:5}

In this section, we list various open questions. The list we present is by no means exhaustive.

\subsection{Random dynamical systems}

In order to model the effect of noise on a discrete-time dynamical system, it is natural to introduce
models obtained by compositions of different maps rather than by repeated applications
of exactly the same transformation. The idea is to study sequences of maps `picked at
random' in some stationary fashion. We refer to \cite[Chap. 5]{handbookb} for a survey.

The simplest case is the following. We assume that the phase space is contained in $\mathds{R}^d$ and that there is a sequence of i.i.d., $\mathds{R}^d$-valued, random variables $\xi_0,\xi_1,\ldots$ such that, instead of observing the orbit of the initial condition $x$, one observes sequences $\{x_n\}$ of points in the state space given by
\[
x_{n+1}=T(x_n)+\epsilon \xi_n
\]
where $\epsilon$ is a fixed parameter (the amplitude of the noise if $|\xi_n|$ is of order one).
The process $\{x_n\}$ is called a {\em stochastic perturbation} of the dynamical system $T$. 
By construction, it is a one-parameter family of Markov chains. If we assume that 
$\xi_n$ has a density $\rho$  with respect to Lebesgue measure, the transition probability
of the chain is given by
\[
p(x_{n+1}| x_n)=\frac{1}{\epsilon}\, \rho\left( \frac{x_{n+1}-T(x_n)}{\epsilon}\right).
\]
One expects that in the limit $\epsilon\to 0$ (the zero-noise limit), the right-hand side converges
to $\delta(x_{n+1}-T(x_n))$ and that, if $\mu_\epsilon$ is an invariant measure for the chain,
then its accumulation points  (in vague topology) should be invariant measures for $T$.
There are reasons to believe that under fairly general conditions, SRB measures may be natural candidates for zero-noise limits, hence they should be stochastically stable. This is indeed proved for Axiom A systems and certain non-uniformly hyperbolic systems, see {\em e.g.} \cite{zeronoise} and \cite{BV} for the H\'enon map. 

A natural question is to prove concentration inequalities for random dynamical systems, in particular
for the additive noise model. This would lead, for instance, to quantitative informations on the
distance between the empirical measure of the process $\{x_n\}$  and the SRB measure $\mu$  
as a function of $n$ and $\epsilon$.

The above setting concerns `dynamical noise'. Another relevant situation is `observational noise':
one observes the process $y_n=x_n+\epsilon \xi_n$ and the goal is merely to extract $\{x_n\}$,
and eventually try to reconstruct $T$ \cite{lalley}.

\subsection{Coupled map lattices}

Coupled map lattices are a class of (discrete-time) spatially extended dynamical systems 
which were introduced in the 1980's by physicists. We refer to the lecture notes \cite{lnp} for more details and
background.

The basic set-up is a state space $\Omega=I^{\mathds{Z}^d}$ where $I\subset\mathds{R}$ is a compact 
interval, typically $[0,1]$.There is a `local' dynamics $\tau:I\circlearrowleft$ which defines an
`unperturbed' dynamics $T_0$ on $\Omega$ by $\big(T_0(x)\big)_i=\tau(x_i)$, $i\in\mathds{Z}^d$. Then
one defines a perturbed dynamics by introducing couplings $\Phi_\epsilon:\Omega\circlearrowleft$
of the form $\Phi_\varepsilon(x)=x+A_\varepsilon(x)$. The basic (and most studied) example is
the `diffusive' nearest neighbor coupling
\[
\big( \Phi_\varepsilon(x)\big)_i=x_i + \frac{\varepsilon}{2d} \sum_{|i-j|=1}(x_j-x_i),\quad i\in\mathds{Z}^d.
\]
Of course, $\varepsilon$ measures the strength of the coupling.\newline
The dynamics we are interested in is 
\[
T_\varepsilon:=\Phi_\epsilon\circ T_0.
\]
The study of such dynamical systems offer many challenges and a lot of questions remain
open \cite{lnp}. 

From the point of view of probabilistic properties, the following is known, see \cite{BGK} and
references therein.
The local map $\tau$ on the unit interval $I$ is assumed to be continuous and piecewise $C^2$. The expansion rate is assumed to be bigger than $2$: $|\tau'|>2$ and both the first- and second-order derivatives are bounded. The couplings are assumed to be diffusive and of finite range (the above example corresponds
to a range equal to one). 
Under these conditions, the coupled map lattice $T_\varepsilon$ has a unique observable measure
$\mu_\epsilon$ in the sense that, for $m^{\otimes \mathds{Z}^d}$-almost every point $x\in\Omega$ state,
\[
\frac{1}{n} \sum_{k=0}^{n-1} \gdelta_{T_\varepsilon^k x} \xrightarrow{\scriptscriptstyle{\textup{vaguely}}}
\mu_{\varepsilon}.
\]
This measure is exponentially mixing both in time and space.
Moreover, any Lipschitz function on $I^{\mathds{Z}^d}$ depending on a finite number of coordinates 
satisfies the central limit theorem with respect to $(T_\epsilon,\mu_\epsilon)$. The authors also prove a local limit theorem. All these results hold provided that $\varepsilon$ is small enough.
As the authors point out, their tools also allow to prove exponential large deviations.

A natural question is to prove concentration inequalities in this context. One expects an exponential
concentration inequality to hold. 

\subsection{Partially hyperbolic systems}

As mentioned above, the theory of hyperbolic dynamical systems initially
developed from the notion of uniform hyperbolicity. This notion can be weakened in essentially
two ways. One of these is to retain hyperbolicity without uniformity, which leads to the theory of non-uniformly hyperbolic dynamical systems. The class of systems modeled by Young towers described in this chapter
is an important subclass of such systems.\newline
The other generalization is to retain uniformity without hyperbolicity by allowing a center direction in which any expansion or contraction is in a uniform way slower than the expansion and contraction in the unstable and stable subspaces. Such systems are called {\em partially hyperbolic}. Among the basic examples
are time-one maps of Anosov flows (the center direction is the flow direction), quasi-hyperbolic toral automorphisms and mostly contracting diffeomorphisms. We refer to \cite[Chap. 1]{handbookb} for a survey.

In  \cite{dolgo}, the author proves many probabilistic results such as the central limit theorem (and its
refinements like the almost-sure invariance principle) and exponential large deviations.\newline
It would be nice to establish concentration inequalities for partially hyperbolic systems.

\subsection{Nonconventional ergodic averages}

Nonconventional or mutiple ergodic averages are typically of the form
\[
\frac{1}{n} \sum_{k=0}^{n-1}
f_1(T^k x)f_2(T^{2k}x)\cdots f_\ell(T^{\ell k}x).
\]
That is, one considers the averages of products of, say, bounded measurable functions
along an arithmetic progression of length $\ell$ for an arbitrary integer $\ell\geq 1$.
The case $\ell=1$ is of course the standard case.
Such averages originated in the ergodic theoretic proof by Furstenberg of Szemer\'edi's theorem
on arithmetic progressions based on the so-called multiple recurrence theorem \cite{furstenberg}.
For a dynamical system $(X,T,\mu)$ which is weakly mixing, the above averages converge
in $L^2$ to $\prod_{k=1}^\ell \int f_k \textup{d}\mu$. \newline
The next questions are about fluctuations of nonconventional averages when the $f_j$'s are, say, Lipschitz
functions : central limit theorem, large
deviations and concentration properties. Regarding the central limit theorem, a first step was
done by Kifer \cite{kifer} for uniformly hyperbolic systems (for averages along more general progressions). 
Large deviations seem much more difficult to analyse and turn out to be nontrivial even for i.i.d. processes
(see \cite{ccgr}).

A transfer operator approach remains to be introduced to tackle such problems because the usual
machinery does not seem appropriate. Remarkably, concentration inequalities, if available for the
system at hand, apply straightforwardly and provide nontrivial informations while they `ignore' the
fine structure of theses averages. We leave as an exercise to the reader
the derivation of such concentration bounds.

\subsection{Erd\"os-R\'enyi law for nonuniformly hyperbolic systems and applications to multifractal analysis}

We come back to large deviations (see Subsection \ref{subsec:LD}). When a rate function does exist for a dynamical system
(see Theorem \ref{LDYexp}), the following question is natural: 
\begin{itemize}
\item[]given an observable $f$, is it possible to extract the rate function $\mathbf{I}_f$ solely from 
a typical orbit of the system ?
\end{itemize}
With a different motivation, this question was answered by Erd\"os and R\'enyi \cite{erdos}
in the context of i.i.d. random variables. In the context of dynamical systems, the first
result was obtained in \cite{chazottes-collet} for a class of piecewise, uniformly expanding
maps of the interval. For this class, Theorem \ref{LDYexp} is valid and one can in fact
get refined large deviation estimates necessary to obtain the following result.
Given an observable $f$ and $t$ in the domain of $\mathbf{I}_f$, let 
\[
M_k(x)=\max\big\{ S_k f(T^jx): 0\leq j\leq\lfloor \exp(k \mathbf{I}_f(t))\rfloor-k\big\}
\]
In words, we are looking for the largest ergodic sum of $f$ in a window
of width $k$ inside the orbit of $x$ up to time $\lfloor \exp(k \mathbf{I}_f(t))\rfloor-k$.

\newpage

\begin{svgraybox}
\begin{theorem}[Erd\"os-R\'enyi law for uniformly expanding maps of the interval \cite{chazottes-collet}]
\label{ERlaw}
\leavevmode\\
Let $T:[0,1]\circlearrowleft$ be  a piecewise $C^2$, uniformly expanding map
which is topologically mixing and $\mu$ its unique absolutely continuous invariant
measure. Let $f:[0,1]\to\mathds{R}$ be an observable of bounded variation\footnote{which is
not of the form $g-g\circ T$ for some bounded measurable $g$.}. Then, there exists $t^*>0$
such that, for any $|t|\leq t^*$ and for Lebesgue-almost every $x\in [0,1]$
\[
\lim_{k\to\infty} \frac{M_k(x)}{k}=t.
\]
More precisely, one has almost everywhere
\[
\limsup_{k\to\infty} \frac{M_k(x)-kt}{\log k}\leq \frac{1}{2u}
\]
and
\[
\liminf_{k\to\infty} \frac{M_k(x)-kt}{\log k}\geq -\frac{1}{2u},
\]
where $u=\mathbf{I}'_f(t)$ .
\end{theorem}
\end{svgraybox}

\medskip

Notice that this theorem gives an optimal rate of convergence, the same as in the i.i.d.
case obtained by Deheuvels {\em et al.} (see \cite{chazottes-collet}).

In view of Theorem \ref{LDYexp} and the technique used in \cite{chazottes-collet}, one
expects that Theorem \ref{ERlaw} be true for systems modeled by a Young tower
with exponential tails. This was partially showed in \cite{DN}, but only in the one-dimensional
case, and with a non-optimal rate.\newline 
On the side of applications, Theorem \ref{LDYexp} allows to construct an
estimator for $\mathbf{I}_f$. This is particularly relevant to the estimation of multifractal
spectra, see \cite{BT}.

\section{Notes on further results}
\label{sec:6}

We quickly describe or barely mention other results that we could not
develop in the main text.

\subsection{More on the central limit theorem}

It is natural to ask for a speed of convergence in the central limit theorem. 
This type of result is called a \emph{Berry-Esseen theorem}.

For systems modeled by a Young tower with exponential tails,
one has the following. Let $f:\Omega\to\mathds{R}$ be a H\"older continuous
observable. Assume that $\sigma_f>0$. Then there exists a constant $c=c(f)>0$ such that
\[
\sup_{t\in\mathds{R}}\left|
\mu\left\{x : \frac{S_nf(x)-n\int f\textup{d}\mu}{\sqrt{n}}\leq t \right\}-
\frac{1}{\sqrt{2\pi}\sigma_f}\int_{-\infty}^t e^{-\frac{u^2}{2\sigma^2_f}} \textup{d}u\right|
\leq \frac{c}{\sqrt{n}},\quad\forall n\in\mathds{N}.
\]
The speed of convergence can be slower. Let us again illustrate this by looking
at the map $T_\alpha$ given by \eqref{MPmap}. For $0<\alpha<1/2$ and $f$
H\"older continuous (which is not of the form $g-g\circ T_\alpha$), we know that
the central limit theorem holds (see end of Section \ref{sec:clt}).\newline
\begin{itemize}
\item If $0<\alpha<1/3$ then one gets a speed of order $\mathcal{O}(1/\sqrt{n})$
as above.
\item If $1/3<\alpha<1/2$ and $f(0)\neq 0$, the speed is
$\mathcal{O}(1/n^{\frac{1}{2\alpha}-1})$. 
\end{itemize}
We refer the interested reader to \cite{seb-BE} for more details and proofs,
where a `local limit theorem' is also proved.

\subsection{Moderate deviations}

One can also characterize the fluctuations of $S_n f$ which are of an order intermediate
between $\sqrt{n}$ (central limit theorem) and $n$ (large deviations). Such fluctuations, when
suitably scaled, satisfy large deviations type estimates with a quadratic rate function
determined by $\sigma_f^2$. We have the following theorem:
\begin{svgraybox}
\begin{theorem}[Moderate deviations \cite{youngLD}]
\leavevmode\\
Let $T: \Omega\circlearrowleft$ be a dynamical system modeled by a
Young tower and $\mu$ its SRB measure.
Assume that  $m^u\{R > n\}= \mathcal{O}(e^{-an})$ for some $a>0$.
Let $f:\Omega\to\mathds{R}$ be a H\"older continuous observable which is not of
the form $g-g\circ T$ (whence $\sigma_f^2>0$).
Let $a_n$ be an increasing sequence of positive real numbers such 
that $\lim_{n\to\infty} a_n/\sqrt{n}=\infty$ and $\lim_{n\to\infty}a_n/n=0$.
Then for any interval $[a,b]\subset \mathds{R}$ we have
\[
\lim_{n\to\infty} \frac{1}{a_n^2/n} \log\mu\left\{x\in\Omega: \frac{S_n f(x)-n\int f\textup{d}\mu}{a_n}\in[a,b] \right\}=-\inf_{t\in[a,b]}\frac{t^2}{2\sigma_f^2}\cdot
\]
\end{theorem}
\end{svgraybox}

\medskip

For the case of systems modeled by Young towers with polynomial tails, see
\cite{melbournePAMS}.

\subsection{Far beyond the CLT: the invariance principle}

The almost sure invariance principle is a very strong reinforcement of the central limit theorem: it ensures that the trajectories of a process can be matched with the trajectories of a Brownian motion in such a way that almost surely the error between the trajectories is negligible compared to the size of the trajectory.

For $\lambda\in (0,1/2]$ and $\Sigma^2$ a (possibly degenerate) symmetric semi-positive-definite $d\times d$ matrix, we say that an $\mathds{R}^d$-valued process $(A_0,A_1,\ldots)$ satisfies an almost sure invariance principle
with error exponent $\lambda$ and limiting covariance $\Sigma^2$ if there exist a probability space $\mathcal{P}$
and two processes $(A^*_0,A^*_1,\ldots)$ and $(B_0,B_1,\ldots)$ on $\mathcal{P}$ such that:
\begin{enumerate}
\item the processes $(A_0,A_1,\ldots)$ and $(A^*_0,A^*_1,\ldots)$ have the same distribution;
\item the random variables $(B_0,B_1,\ldots)$ are independent and distributed as $\mathcal{N}_{0,\Sigma^2}$;
\item and almost surely in $\mathcal{P}$
\[
\left|\sum_{\ell=0}^{n-1} A^*_\ell -\sum_{\ell=0}^{n-1} B_\ell\right| = o(n^\lambda).
\]
\end{enumerate}
A Brownian motion at integer times coincides with a sum of i.i.d. Gaussian variables, hence this definition can also be formulated as an almost sure approximation by a Brownian motion, with error $o(n^\lambda)$.

In the dynamical system context, take $A_\ell=f\circ T^\ell$ where $f:\Omega\to\mathds{R}^d$ is
regular. It is proved  in \cite{melbourneasip} by martingale methods and then in  \cite{gouezelasip}
with purely spectral methods, that a dynamical systems modeled by Young towers satisfy
the almost-sure invariance principle. Namely, this is the case if $\int R^q \textup{d}m^u<\infty$
for $q>2$ and for observables $f:\Omega\to\mathds{R}^d$ which are H\"older continuous.
The relevance of considering $\mathds{R}^d$-valued observable is that, for instance,
the position variable of the planar periodic Lorentz gas with finite horizon approximates a two-dimensional Brownian motion.

The almost-sure invariance principle implies in particular the central limit theorem, the functional
central limit theorem, and the law of
iterated logarithm, among others, see {\em e.g.} \cite{hall,PS}. It also
implies the almost-sure central limit theorem \cite{lacey}.
 
\begin{acknowledgement}
The author thanks S\'ebastien Gou\"ezel for useful comments. He also thanks
Cesar Maldonado and Mike Todd for a careful reading.
\end{acknowledgement}

\end{document}